\documentclass{conm-p-l}
\usepackage[all, cmtip]{xy}
\usepackage{amsthm}
\usepackage{xcolor} 
\usepackage{mathrsfs} 
\usepackage[T1]{fontenc}

\newtheorem{theorem}{Theorem}[section]

\newtheorem{proposition}[theorem]{Proposition}
\newtheorem*{recipe}{Secondarization Recipe}

\theoremstyle{definition}
\newtheorem{definition}[theorem]{Definition}
\newtheorem{example}[theorem]{Example}

\theoremstyle{remark}
\newtheorem{remark}[theorem]{Remark}

\numberwithin{equation}{section}

\DeclareFontFamily{OT1}{wncyi}{}
\DeclareFontShape{OT1}{wncyi}{m}{it}{
<5> <6> <7> <8> <9> gen * wncyi
<10> <10.95> <12> <14.4> <17.28> <20.74> <24.88> wncyi10
}{}
\DeclareSymbolFont{cyrletters}{OT1}{wncyi}{m}{it}
\DeclareSymbolFontAlphabet{\cyrmath}{cyrletters}
\DeclareMathSymbol{\rE}{\cyrmath}{cyrletters}{003}

\makeatletter 
\def\l@subsection{\@tocline{2}{0pt}{1pc}{5pc}{}} \def\l@subsection{\@tocline{2}{0pt}{2pc}{6pc}{}}
\makeatother

\begin{document}

\title[Vinogradov's Cohomological Geometry of PDEs]{Vinogradov's Cohomological Geometry\\of Partial Differential Equations}

\author{Fabrizio Pugliese}
\address{DipMat, Universit\`a degli Studi di Salerno, via Giovanni Paolo II n${}^{\circ}$123, 84084 Fisciano (SA), Italy.}
%\curraddr{Department of Mathematics and Statistics,
%Case Western Reserve University, Cleveland, Ohio 43403}
\email{fpugliese@unisa.it}

\author{Giovanni Sparano}
\address{DipMat, Universit\`a degli Studi di Salerno, via Giovanni Paolo II n${}^{\circ}$123, 84084 Fisciano (SA), Italy.}
\email{sparano@unisa.it}

\author{Luca Vitagliano}
\address{DipMat, Universit\`a degli Studi di Salerno, via Giovanni Paolo II n${}^{\circ}$123, 84084 Fisciano (SA), Italy.}
\email{lvitagliano@unisa.it}

%    General info
\date{}

\dedicatory{to the memory of our teacher and colleague Alexandre Mikhailovich Vinogradov}


\begin{abstract}
Secondary Calculus is a formal replacement for differential calculus on the space of solutions of a system of possibly non-linear partial differential equations and it is essentially due to Alexandre M. Vinogradov and his collaborators. Many coordinate free properties of PDEs find their natural place in Secondary Calculus including: symmetries and conservation laws, variational principles and the coordinate free aspects of the calculus of variations, recursion operators and Hamiltonian structures, etc. The building blocks of this language are horizontal cohomologies of diffieties, i.e. infinite prolongations of PDEs, and their versions with local coefficients. The main paradigm of Secondary Calculus is the principle, due to A. M. Vinogradov, roughly stating that: differential calculus on the space of solutions of a PDE is calculus up to homotopy on the horizontal De Rham algebra of the associated diffiety. We will review the fundamentals of Secondary Calculus including its main motivations. In the last part of the paper, we will try to explain the role of homotopy in the theory.
\end{abstract}

\maketitle

\section*{Introduction}

Systems of algebraic equations can be encoded geometrically in an algebraic variety. This apparently harmless observation is the starting point for the development of such a successful branch of modern Mathematics as Algebraic Geometry. One of the leading principles behind most of the mathematical work of Alexandre M. Vinogradov is that Partial Differential Equations (PDEs) should be treated in a similar way. Namely, let $x = (x^1, \ldots, x^n)$ be some independent variables, let $u = (u^1, \ldots, u^m)$ be some dependent variables, and let 
\begin{equation}\label{eq:F=0}
\mathscr E_0 : F_a (x, \ldots, u_I, \ldots) = 0,
\end{equation}
$a = 1, \ldots, r$, be a system of possibly non-linear PDEs, where a multi-index $I$ denotes multiple partial derivatives with respect to the $x$'s. Adding to (\ref{eq:F=0}) all its total derivatives results in a new but equivalent system of infinitely many PDEs
\begin{equation}\label{eq:DF=0}
D_J F_a (x, \ldots, u_I, \ldots) = 0.
\end{equation}
Here $D_J$ denotes multiple total derivatives of arbitrarily high order with respect to the $x$'s. Now, interpret $(x, \ldots, u_I, \ldots)$ as independent coordinates on an $\infty$-dimensional manifold $J^\infty$. Then (\ref{eq:DF=0}) can be interpreted as defining a (generically $\infty$-dimensional) submanifold $\mathscr E \subseteq J^\infty$, the $\infty$-prolongation of $\mathscr E_0$. This submanifold is canonically equipped with an involutive distribution $\mathscr C$, called the \emph{Cartan distribution}, spanned by the total derivatives $D_i|_{\mathscr E}$, $i =1, \ldots, n$. The pair $(\mathscr E, \mathscr C)$ is, in the terminology of Vinogradov and his school, a \emph{diffiety}\footnote{To the best of our knowledge, the term \emph{diffiety} appeared for the first time in \cite{Vin1984}.} (from the words \emph{diff}erential equation and var\emph{iety}) and contains most of the relevant information about the original system $\mathscr E_0$ (see, e.g., the extensive monographs \cite{b...99, kv98, v01} on the subject). According to Vinogradov, the diffiety $(\mathscr E, \mathscr C)$ should be considered as the ultimate geometric portrait of $\mathscr E_0$. In other words, diffieties are to PDEs what algebraic varieties are to algebraic equations. However the geometry of a diffiety is much more intricate than that of an algebraic variety. Namely, while solutions of a system of algebraic equations are points of an algebraic variety and they do not possess any internal structure, solutions of a system of PDEs $\mathscr E_0$ are integral submanifolds of the Cartan distribution on the diffiety $(\mathscr E, \mathscr C)$ and they do possess internal structure. As a consequence the space of solutions of a PDE share some features of the leaf space of a foliation (and \emph{it is} the leaf space of a foliation when $\dim \mathscr E < \infty$) and should be treated as a stack rather than as an ordinary space \cite{mm03, BX2011, dH2013}. This is exactly where homological and homotopical algebra come into play when trying to develop a general theory of PDEs. A fundamental achievement of Vinogradov is that there is a cohomological replacement for \emph{differential calculus on the space of solutions of a PDEs}. Vinogradov used to call such formal calculus \emph{Secondary Calculus}, because of an analogy with secondary quantization in physical field theory \cite{v98, b...99, v01}. The main building blocks of Secondary Calculus are the leaf-wise cohomologies of the Cartan distribution, often called \emph{horizontal De Rham cohomologies}. Most of the coordinate free properties of a PDE can be expressed in terms of horizontal cohomologies: symmetries, conservation laws, variational principles, etc., to name a few. Additionally, horizontal cohomologies can be interpreted as smooth functions, vector fields, differential forms on the space of solutions of the PDE. One of the main features of horizontal cohomologies supporting this interpretation is that they possess the \emph{correct} algebraic structures (those that we expect from functions, vector fields, differential forms). In summary, Secondary Calculus is the collection of all these algebraic structures. 

In the ultimate idea of Vinogradov this fascinating constructions should be studied in their homotopy aspects. More precisely, Vinogradov conjectured that what is important is not really horizontal cohomology but rather the homotopy type of the horizontal De Rham algebra, and the homotopy category (or even the derived category) of differential graded modules over it \cite[Chapter 5]{v01}. Notice that this conjecture (together with the interpretation of the space of solutions of a PDE as a stack), is somehow complementary to the interpretation of a PDE as a derived zero locus and, in the particular case of an Euler-Lagrange PDE, as a derived critical locus (see, e.g., \cite{v20}, and references therein), and suggests a relationship to Homotopical and Derived Geometry \cite{tv05, tv08, t14} that, unfortunately, to the best of our knowledge, has not been yet investigated. At this stage, we can only speculate that the space of solutions of a PDE should be ultimately treated as a \emph{higher derived stack} \cite{t10}.

In this exposition we will focus on the theory rather than on applications or the explicit computation of the involved cohomologies. We will omit the proofs. They can be found in the cited references. The paper is divided into three sections: In Section \ref{Sec:Geom} (Geometry of PDEs) we explain how to encode a system of PDEs in a geometric object: a \emph{diffiety}. We also discuss the main geometric structure possessed by a diffiety, namely its Cartan distribution. Our main sources for this section are \cite{b...99, kv98, klv86}. In Section \ref{Sec:Homol} (Homology of PDEs) we show that there is a natural cohomology attached to a diffiety (hence to a PDE), the \emph{horizontal cohomology}. This is the core of the paper, where Vinogradov's Secondary Calculus is discussed. Horizontal cohomology is a generalization of the leaf-wise cohomology of a foliated manifold and contains important coordinate free information on a PDE. We provide an interpretation of the main horizontal cohomologies. Horizontal cohomologies can also be seen as the building blocks of a \emph{differential calculus on the space of solutions of a PDE}, what we call \emph{Secondary Calculus}. One of the main supporting facts for the latter interpretation is that horizontal cohomologies are equipped with the \emph{correct} algebraic structures for smooth functions, vector fields, differential forms, etc. Our main sources here are \cite{b...99, kv98, v01}, see also \cite{v84, v09}. In Section \ref{Sec:Homot} (Homotopy of PDEs) we present the latest developments. This section is an extremely compact review of the papers \cite{vi12, v12} by the third author. We show that the main algebraic structures on horizontal cohomologies do all come from appropriate \emph{structures up to homotopy} on horizontal cochains. This supports Vinogradov's idea that the appropriate category to develop Secondary Calculus is (a subcategory of) the derived category of DG modules over the horizontal De Rham algebra.

\tableofcontents

\section{Geometry of PDEs}\label{Sec:Geom}

In this section we show that a system of, generically non-linear, PDEs can be encoded into a geometric object: a \emph{diffiety}. We will work in the rather general setting of a PDEs imposed on submanifolds $N$ of a fixed dimension $n$ in a manifold $P$ (see \cite{klv86}), in contrast with most of the literature where only PDEs imposed on sections of a fibration (or a fiber bundle) are considered.

\subsection{jet spaces}
Let $n,m$ be non-negative integers, and let $P$ be a manifold of dimension $n + m$. We begin providing an intrinsic definition of \emph{derivatives} of an $n$-dimensional submanifold $N \subseteq P$. In the following, $n$ should be interpreted as the number of \emph{independent variables}, $m$ as the number of \emph{dependent variables}, and $P$ as the space parameterizing both independent and dependent variables. In other words, in this theory, dependent and independent variables can be mixed and swapped. Around every point $p \in N$ there are coordinates $(x^i, u^\alpha)$ on $P$, $i = 1, \ldots, n$, $\alpha = 1,\ldots, m$, \emph{adapted} to $N$ in the sense that, in these coordinates, $N$ looks like the \emph{graph of a map}:
\begin{equation}\label{eq:N_loc}
N :  u^\alpha = f^\alpha (x)  .
\end{equation}
We will need to take multiple partial derivatives of the $f$'s with respect to the $x$'s. We adopt the following notation. Let $h$ be a non-negative integer and let $I = i_1\cdots i_h$ be a \emph{multiindex}, i.e.{} a (possibly empty) word containing $h$ letters $i_\ell \in \{1, \ldots, n\}$, $\ell = 1, \ldots, h$. We identify two words if they only differ by a permutation of their letters. We also compose two words by concatenation. Denote $|I| := h$ and call it the \emph{lenght} of the multiindex $I$. Finally denote
\[
\frac{\partial^{|I|} f^\alpha}{\partial x^I} := \frac{\partial^h f^\alpha}{\partial x^{i_1} \ldots \partial x^{i_h}}.
\] 

If there is another $n$-dimensional submanifolds $\tilde N$ through $p$, i.e.{} $p \in N \cap \tilde N$, then around $p$ there are coordinates $(x^i, u^\alpha)$ on $P$ adapted to both $N, \tilde N$:
\[
N : u^\alpha = f^\alpha (x), \quad \text{and} \quad \tilde N : u^\alpha = \tilde f{}^\alpha (x) .
\]
In this case we say that $N$ and $\tilde N$ are \emph{tangent up to order $k = 0, 1, \ldots, \infty$ at $p$} if
\[
\frac{\partial^{|I|} f^\alpha}{\partial x^I} (p) = \frac{\partial^{|I|} \tilde f{}^\alpha}{\partial x^I} (p),\quad \text{for all multiindexes $I$ with $|I| \leq k$},
\]
and we write $N \sim_p^k \tilde N$. Tangency up to order $k$ at $p$ is a well-defined equivalence relation on submanifolds through $p$. In particular, it is independent of the choice of adapted coordinates. The equivalence class of $N$ is denoted $N^k_p$ and called the \emph{$k$-jet of $N$ at $p$}.   In other words,   $k$-jets at $p$ encode partial derivatives of $n$-dimensional submanifolds at $p$ up to order $k$. The space of tangency classes is denoted $J^k_p (P, n)$. Finally put
\[
J^k (P, n) = \coprod_{p \in P} J^k_p (P,n).
\] 
We call $J^k (P, n)$ the \emph{$k$-jet space of $n$-dimensional submanifolds of $P$}. It is clear that $J^0 (P, n)$ identifies naturally with $P$ and $J^1 (P, n)$ identify with the bundle of $n$-Grassmannians in the fibers of the tangent bundle $TP$. More precisely, the $0$-jet of $N$ at $p$ identifies with $p$ itself, while the $1$-jet of $N$ at $p$ identifies with the tangent space $T_p N$. For a general $k$, the space $J^k (P, n)$ can be coordinatized as follows. Choose coordinates $(x^i, u^\alpha)$ on $P$, and let $U \subseteq J^k (P, n)$ be the subset consisting of $k$-jets of submanifolds for which $(x^i, u^\alpha)$ are adapted coordinates. We define coordinates  $(x^i, u^\alpha_I)$ on $U$ by putting
\[
x^i \big( N^k_p \big) := x^i (p), \quad \text{and} \quad u^\alpha_I \big(N^k_p\big) = \frac{\partial^{|I|} f^\alpha}{\partial x^I} (p),
\]
for all $N$ locally given by (\ref{eq:N_loc}), $i = 1, \ldots, n$, $\alpha = 1, \ldots, m$, and $|I| \leq k$. The $(x^i, u^\alpha)$ are called \emph{jet coordinates}, they cover $J^k (P, n)$ and, for $k < \infty$, they give to $J^k (P, n)$ the structure of a smooth manifold. For all $h \leq k \leq \infty$ the are projections
\[
J^k (P, n) \to J^h (P, n), \quad N^k_p \mapsto N^h_p,
\]
and by Borel Lemma $J^\infty (P, n)$ identifies with the inverse limit of the tower of surjective submersions
\[
J^0 (P, n)  \longleftarrow \cdots \longleftarrow J^{k-1} (P, n) \longleftarrow J^k (P,n) \longleftarrow \cdots.
\]
%$J^\infty (P, n) \cong \underset{\longleftarrow}{\lim} J^k (P, n)$.
This gives to the $\infty$-jet space the structure of a pro-finite dimensional manifold (see, e.g., \cite{b...99,v01}, see also \cite{g17}).

\begin{remark}
For the reader more familiar with jets of sections, we mention that, when $P$ is fibered over some $n$-dimensional manifold $M$, then, for all $k \leq \infty$, the space $J^k (P, M)$ of $k$-jets of sections of $P$ over $M$ can be recovered as the open and dense submanifold of $J^k (P, n)$ consisting of jets of $n$-dimensional submanifolds transverse to the projection $P \to M$. The open embedding $J^k (P, M) \hookrightarrow J^k (P, n)$ maps the $k$-jet at $x \in M$ of a (possibly local) section $s$ of $P \to M$ to the $k$-jet of its image at $s(x)$.
\end{remark}

Differential calculus on $J^\infty (P, n)$ can be defined algebraically. For instance, we define smooth functions $C^\infty (J^\infty (P, n))$ on $J^\infty (P,n)$ as the direct limit of the algebra filtration
\begin{equation}\label{eq:filt_C}
C^\infty \big(J^0 (P, n)\big) \hookrightarrow \cdots \hookrightarrow C^\infty \big(J^{k-1} (P, n)\big) \hookrightarrow C^\infty \big(J^k (P, n)\big) \hookrightarrow \cdots.
\end{equation}
In other words a smooth function on $J^\infty (P, n)$ is a smooth function of the independent variables $x^i$, the dependent variables $u^\alpha$, and their derivatives up to some finite order $k$, but $k$ can be arbitrarily high (sometimes they are defined as functions which are only locally of this type). Vector fields, differential forms, etc., on $J^\infty (P, n)$ are defined in a way compatible with the filtration (\ref{eq:filt_C}). We will not insist on these technicalities and we refer to \cite{b...99} for details.

\begin{remark}\label{rem:infty_dim}
We stress that, in this $\infty$-dimensional setting, some important results in finite dimensional differential geometry might fail. To mention a few: the inverse function theorem, existence and uniqueness of the flow of a vector field, the Frobenius theorem, etc.
\end{remark}

Any $n$-dimensional submanifold $N \subseteq P$ defines a new $n$-dimensional submanifold
\[
N^k :=\left\{ N^k_p : p \in N \right\}
\]
in $J^k (P, n)$, called the \emph{$k$-jet prolongation} of $N$. If $N$ locally looks like (\ref{eq:N_loc}) in coordinates, then $N^k$ locally looks like
\[
N^k :  u^\alpha_I = \frac{\partial^{|I|} f^\alpha}{\partial x^I} (x), \quad |I| \leq k,
\]
in jet coordinates.

\begin{remark}\label{rem:emb_N}
Let $N \subseteq P$ be an $n$-dimensional submanifold. For all $k \leq \infty$, the map
\[
j^k : N \to N^k, \quad p \mapsto N^k_p
\]
is a diffeomorphism. In particular $N$ embeds into $J^k (P, n)$ for all $k$.  
When $P$ is fibered over some $n$-dimensional manifold $M$, and $s : M \to P$ is a (possibly local) section of $P$ over $N$, then the usual $k$-jet prolongation $j^k s$ of $s$ can be recovered as the composition
\[
j^k s : M \overset{s}{\longrightarrow} \operatorname{im} s \overset{j^k}{\longrightarrow} (\operatorname{im} s)^k \hookrightarrow J^k (P, M).
\]
  
\end{remark}

\begin{remark}\label{rem:jet_emb}
Let $r$ be a non-negative integer and let $Q \subseteq P$ be an $(n+r)$-dimensional submanifold. Then, for all $k \geq 0$, the $k$-jet space $J^k (Q, n)$ is a submanifold in the $k$-jet space $J^k (P, n)$ in the obvious natural way.
\end{remark}

\subsection{the Cartan distribution}\label{Sec:CD}
The main geometric structure on jet spaces is the \emph{Cartan distribution}. There are several equivalent ways to define it. The one we present here will be also useful in defining the prolongation of a PDE in the next subsection. We begin noticing that a section of the fibration $J^1 (P, n) \to J^0 (P, n) \cong P$ can be interpreted as an $n$-dimensional distribution on $P$. Now, for every $1 \leq k \leq \infty$, the $k$-jet space is canonically equipped with a rank $n$ \emph{distribution along the projection $J^{k} (P, n) \to J^{k-1}(P, n)$}, i.e.{} a smooth map
\[
\mathscr C : J^k (P, n) \longrightarrow J^1 \big(J^{k-1} (P, n), n\big)
\]
such that the following diagram commutes
\[
\begin{array}{c}
\xymatrix{
 & J^1 \big(J^{k-1} (P, n), n\big) \ar[d] \\
 J^k (P, n) \ar[ur]^-{\mathscr C} \ar[r] & J^{k-1} (P, n)
}
\end{array}.
\]
The \emph{Cartan distribution} $\mathscr C$ is defined as follows. For $z = N^k_p \in J^k (P, n)$, denote by $\bar z = N^{k-1}_p \in J^{k-1}(P, n)$ its projection. Now put
\[
\mathscr C_z := \big( N^{k-1} \big){}^1_{\bar z} = T_{\bar z} \big(N^{k-1} \big) \in J^1 \big(J^{k-1} (P, n), n\big).
\]
In other words $\mathscr C_z$ is the $1$-jet at $\bar z$ of the $(k-1)$-jet prolongation of $N$. In local coordinates
\begin{equation}\label{eq:C_loc}
\mathscr C^\ast ( u^\alpha_{I, j}) = u^\alpha_{Ij}, \quad |I| \leq k-1, \quad j = 1, \ldots, n,
\end{equation}
where, in the left hand side, we denoted by $(x^i, u^\alpha_{I, j})$ the jet coordinates on the $1$-jet space $J^1 \big(J^{k-1} (P, n), n\big)$ (corresponding to the coordinates $(x^i, u^\alpha_I)$ on $J^{k-1}(P, n)$), and, in the right hand side, we are using concatenation of multiindexes.

\begin{remark}\label{rem:J^k+h}
The Cartan distribution $\mathscr C : J^k (P, n) \longrightarrow J^1 \big(J^{k-1} (P, n), n\big)$ is an embedding that we often use to interpret the $k$-jet space $J^k (P, n)$ as a submanifold in the jet space $J^1 \big(J^{k-1} (P, n), n\big)$. More generally for any two non-negative integers $h, k$, there is an embedding
\[
J^{h+k} (P, n) \hookrightarrow J^h \big(J^{k} (P, n), n\big)
\]
given by
\[
N^{h+k}_p \mapsto (N^k)^h_{\overline p}, \quad \overline p = N^k_p.
\]
In the next subsection we will use this embedding to interpret $J^{h+k} (P, n)$ as a submanifold in $J^h \big(J^{k} (P, n), n\big)$.
\end{remark}

For $k = \infty$, the Cartan distribution is a honest rank $n$ distribution on $J^\infty (P, n)$, and (\ref{eq:C_loc}) shows that it is locally spanned by the \emph{total derivatives}, i.e.{} the following vector fields
\[
D_i := \frac{\partial}{\partial x^i} + \sum_{|I| < \infty} u^\alpha_{Ii} \frac{\partial}{\partial u^\alpha_I}, \quad i = 1, \ldots, n.
\]
Dually, the annihilator $\mathscr C^0 \subseteq T^\ast J^\infty (P, n)$ is spanned by the \emph{Cartan forms}
\[
du^\alpha_I - u^\alpha_{Ii}dx^i, \quad \alpha = 1, \ldots, m, \quad |I| \leq \infty.
\]
 
Here, and in the rest of the paper, we adopt the Einstein summation convention over pairs of lowercase-uppercase indexes. We will not adopt the Einstein convention for multiindexes.

\begin{proposition}
The Cartan distribution on $J^\infty (P, n)$ is involutive, i.e.{} the commutator of two sections of $\mathscr C$ lays in $\mathscr C$ as well. Additionally it detects $\infty$-jet prolongations in the following sense: a submanifold $S \subseteq J^\infty (P, n)$ is locally of the form  $S = N^\infty$ for some $n$-dimensional submanifold $N \subseteq P$ if and only if it is an $n$-dimensional integral submanifold of $\mathscr C$. 
\end{proposition}

So, morally, the $\infty$-jet prolongation construction $N \mapsto N^\infty$ identifies $n$-dimensional submanifolds of $P$ with $n$-dimensional integral submanifolds of $\mathscr C$. In the following, \emph{Cartan distribution} will always mean the Cartan distribution on $J^\infty (P, n)$ (unless otherwise stated).

For their importance in the next section, we now discuss \emph{infinitesimal symmetries of the Cartan distribution $\mathscr C$}, i.e.{} vector fields $Y$ on $J^\infty (P, n)$ such that $[Y, \Gamma (\mathscr C)] \subseteq \Gamma (\mathscr C)$. Denote by $\mathfrak X_{\mathscr C}$ the Lie algebra of infinitesimal symmetries of $\mathscr C$ and notice that, by involutivity, $\Gamma (\mathscr C) \subseteq \mathfrak X_{\mathscr C}$ is an ideal. We want to show that the quotient Lie algebra $\mathfrak X_{\mathscr C} / \Gamma (\mathscr C)$ identifies with sections of an appropriate vector bundle $V \to J^\infty (P, n)$. We first define $V$. So let $z = N^\infty_p \in J^\infty (P, n)$. The fiber of $V$ over $z$ is then the $m$-dimensional vector space
\[
V_z := T_p P / T_p N.
\]
In the following, we denote by $\pi : J^\infty (P, n) \to P$ the projection. Now, let $Y \in \mathfrak X_{\mathscr C}$. We define a section $\overline Y$ of $V$ as follows. For $z = N^\infty_p \in J^\infty (P, n)$ put
\[
\overline Y_z := d\pi (Y_z) \operatorname{mod} T_p N \in T_p P / T_p N = V_z.
\]

\begin{proposition}
The assignment $Y \mapsto \overline Y$ defines an $\mathbb R$-linear surjection $\mathfrak X_{\mathscr{C}} \to \Gamma (V)$ with kernel given by $\Gamma (\mathscr C)$. Hence there is a canonical short exact sequence of vector spaces
\[
0 \longrightarrow \Gamma (\mathscr C) \longrightarrow \mathfrak X_{\mathscr C} \longrightarrow \Gamma (V) \longrightarrow 0.
\]
In particular, $\Gamma (V)$ inherits from $\mathfrak X_{\mathscr C}$ a Lie bracket.
\end{proposition}

The Lie bracket on $\Gamma (V)$ is sometimes called the \emph{higher Jacobi bracket} and denoted by $\{-,-\}$. We conclude this section describing it in local coordinates. First notice that $\Gamma (V)$ is locally spanned by the sections $v_\alpha$ defined as follows. For $z = N^\infty_p \in J^\infty (P, n)$ put
\[
v_\alpha |_z := \frac{\partial}{\partial u^\alpha}|_p \operatorname{mod} T_p N \in T_p P / T_p N = V_z.
\]
The isomorphism $\Gamma (V) \cong \mathfrak X_{\mathscr{C}} / \Gamma (\mathscr C)$ is now locally given by
\[
\chi \mapsto \rE_\chi \operatorname{mod} \Gamma (\mathscr C)
\]
where, for $\chi = \chi^\alpha v_\alpha \in \Gamma (V)$, we denoted by $\rE_\chi$ the (local) infinitesimal symmetry of $\mathscr C$ locally given by
\begin{equation}\label{eq:ev_vf}
\rE_\chi = \sum_{|I| < \infty} D_I \chi^\alpha \frac{\partial}{\partial u^\alpha_I},
\end{equation}
and, for any multiindex $i_1\cdots i_h$,
\[
D_{i_1 \cdots i_h} := D_{i_1} \circ \cdots \circ D_{i_h}.
\]
 
We mention for the reader unfamiliar with the cyrillic alphabet, that the script $\rE$ appearing in (\ref{eq:ev_vf}) is pronounced [e], like in s\underline{e}t.
  
A vector field of the form $\rE_\chi$ is sometimes called an \emph{evolutionary vector field}. Notice that evolutionary vector fields are only locally defined on $\infty$-jets of submanifolds (while one can give a coordinate free meaning to (\ref{eq:ev_vf}) on $\infty$-jets of sections of a fibration,  see e.g.{} \cite[Chapter 4, Section 2.4]{b...99} or \cite[Section 3.9]{kv98}  ). It easily follows from (\ref{eq:ev_vf}) that, if $\chi$ and $\psi$ are sections of $V$ locally given by $\chi = \chi^\alpha v_\alpha$ and $\psi = \psi^\beta v_\beta$, then their higher Jacobi bracket $\{ \chi, \psi\}$ is the section locally given by
\[
\{ \chi, \psi \} = \big(\rE_\chi \psi^\alpha - \rE_\psi \chi^\alpha\big) v_\alpha = \sum_{|I|< \infty}\left(D_I \chi^\beta \frac{\partial \psi^\alpha}{\partial u^\beta_I} - D_I \psi^\beta \frac{\partial \chi^\alpha}{\partial u^\beta_I}\right) v_\alpha.
\]

\subsection{diffieties} We now use jets to define PDEs and present their geometric portraits: diffieties. Let $P$ be a manifold of dimension $n + m$.

\begin{definition}\label{def:PDE}
A \emph{system of $k$-th order partial differential equations imposed on $n$-dimensional submanifolds of $P$}, or simply a \emph{PDE}, is a closed submanifold $\mathscr E_0 \subseteq J^k (P, n)$ of the $k$-jet space. A \emph{solution} of a PDE $\mathscr E_0 \subseteq J^k (P, n)$ is an $n$-dimensional submanifold $N \subseteq P$ such that $N^k \subseteq \mathscr E_0$.
\end{definition}

\begin{remark}
As a minimal regularity condition on a PDE $\mathscr E_0 \subseteq J^k (P, n)$ we will always assume $\dim \mathscr E_0 \geq n$ so that the existence of solutions is not excluded a priori by trivial dimensional reasons.
\end{remark}

For a PDE $\mathscr E_0 \subseteq J^k (P, n)$, locally, in jet coordinates, we have
\begin{equation}\label{eq:PDE}
\mathscr E_0 : F_a (x, \ldots, u_I, \ldots) = 0, \quad |I| \leq k,
\end{equation}
for some local smooth functions $F_a \in C^\infty (J^k (P, n))$. If $N \subseteq P$ is an $n$-dimensional submanifold locally given by
\[
N : u^\alpha = f^\alpha (x) ,
\]
then $N$ is a solution of $\mathscr E_0$ iff
\[
F_a \left(x, \ldots, \frac{\partial^{|I|}f}{\partial x^I}, \ldots\right) = 0 .
\]
This motivates Definition \ref{def:PDE}.

In Algebraic Geometry, given a system of algebraic equations, it is natural to take all their algebraic consequences. In other words, given a subset of polynomials, it is natural to consider the ideal that they span. Similarly, given a PDE (\ref{eq:PDE}), it is natural (besides their smooth consequences) to take all their total derivatives, this can be done in a coordinate free way via a construction known as \emph{prolongation} that we now describe. Let $\mathscr E_0 \subseteq J^k (P, n)$ be a PDE, and let $q = 0, 1, \ldots, \infty$.

\begin{definition}\label{def:prolong}
The \emph{$q$-prolongation} of $\mathscr E_0$ is the subset $\mathscr E_q \subseteq J^{k+q} (P, n)$ of the $(k+q)$-jet space defined by
\begin{equation}\label{eq:prolong}
\mathscr E_q := J^q (\mathscr E_0, n) \cap J^{k+q} (P, n).
\end{equation}  
The $\infty$-prolongation will be often denoted simply $\mathscr E$.
\end{definition}

Definition \ref{def:prolong} requires some little explanations. As $\mathscr E_0$ is a submanifold in $J^k (P, n)$, its jet space $J^q (\mathscr E_0, n)$ is a submanifold in the jet space $J^q (J^k (P, n), n)$ as in Remark \ref{rem:jet_emb}. The jet space $J^{k+q} (P, n)$ can be seen as a submanifold in $J^q (J^k (P, n), n)$ as well via the embedding in Remark \ref{rem:J^k+h}. This explains the intersection in (\ref{eq:prolong}). We stress that we will always understand the prolongation $\mathscr E_q$ as a subset in $J^{k+q} (P, n)$. If $\mathscr E_0$ is locally given by (\ref{eq:PDE}), then its $q$-prolongation is locally given by
\[
\mathscr E_q : D_J F_a (x, \ldots, u_I, \ldots) = 0, \quad |J| \leq q.
\]
This shows that $\mathscr E_q$ encodes in a coordinate free way the total derivatives of $\mathscr E_0$ up to order $q$, as desired.

\begin{remark}
For a generic $\mathscr E_0 \subseteq J^k (P, n)$ its $q$-th prolongation $\mathscr E_q$ is \emph{not} a submanifold in $J^{k + h} (P, n)$. When it is a submanifold, then it is clearly a PDE with the same solutions as $\mathscr E_0$. In the following, we will always assume that the $\infty$-th prolongation $\mathscr E = \mathscr E_\infty$ is a submanifold in $J^\infty (P, n)$. This happens, e.g., when $\mathscr E_0$ is a \emph{formally integrable PDE}. Our assumptions are not a great loss of generality according to \emph{Cartan-Kuranishi Prolongation Theorem} which roughly states that, under mild regularity conditions, every PDE becomes formally integrable after finitely many $1$-prolongations.
\end{remark}

Now on we concentrate on the $\infty$-prolongation $\mathscr E \subseteq J^\infty (P, n)$.

\begin{proposition}\label{prop:diffiety}
Let $\mathscr E_0 \subseteq J^k (P, n)$ be a PDE. Assume that its $\infty$-prolongation $\mathscr{E} \subseteq J^\infty (P, n)$ is a submanifold. Then the Cartan distribution restricts to $\mathscr E$ in the sense that
\[
\mathscr C_z \subseteq T_z \mathscr E \quad \text{for all $z \in \mathscr E$}.
\]
Hence the restriction $\mathscr C : \mathscr E \to J^1 (\mathscr E, n)$ is a well-defined rank $n$ involutive distribution on $\mathscr E$. Additionally it detects solutions of $\mathscr E_0$ in the following sense: a submanifold $S \subseteq \mathscr E$ is locally of the form $S = N^\infty$ for some solution $N$ of $\mathscr E_0$ if and only if it is an $n$-dimensional integral submanifold of $\mathscr C$. 
\end{proposition}

The restriction of the Cartan distribution to the $\infty$-prolongation of a PDE will be denoted by $\mathscr C$ as well and called the \emph{Cartan distribution}.

\begin{definition}
A \emph{diffiety} (from the words \emph{differential equation} and \emph{variety}) is a pair $(\mathscr E, \mathscr C)$ where $\mathscr E$ is the $\infty$-prolongation of some PDE and $\mathscr C$ is the Cartan distribution on $\mathscr E$.
\end{definition}

Proposition \ref{prop:diffiety} now shows that a diffiety $(\mathscr E, \mathscr C)$ contains most of the relevant information about the PDE $\mathscr E_0$ defining it.

We conclude this section briefly discussing \emph{infinitesimal symmetries} of a PDE. Let $(\mathscr E, \mathscr C)$ be the diffiety corresponding to a PDE $\mathscr E_0$.

\begin{definition}\label{def:inf_symm_E}
An \emph{infinitesimal symmetry} of $\mathscr E_0$ is a vector field $Y \in \mathfrak X (\mathscr E)$ such that $[Y, \Gamma (\mathscr C)] \subseteq \Gamma (\mathscr C)$. A \emph{non-trivial infinitesimal symmetry} is an infinitesimal symmetry modulo \emph{trivial ones}, i.e.{} sections of $\mathscr C$.
\end{definition}

\begin{remark}
The terminology ``trivial, non-trivial symmetries'' in Definition \ref{def:inf_symm_E} is motivated by the following facts. Let $X$ be a section of the Cartan distribution on $J^\infty (P, n)$ and let $(\mathscr E, \mathscr C)$ be any diffiety, then $X$ is tangent to $\mathscr E$ and $X|_{\mathscr E}$ is an infinitesimal symmetry of $\mathscr E_0$. In other words, sections of the Cartan distribution are infinitesimal symmetries of \emph{all} PDEs (regardless what is the PDE we are considering). 
\end{remark}

There is an efficient description of (non-trivial) infinitesimal symmetries in local coordinates. Consider the vector bundle $V \to J^\infty (P, n)$ from the previous subsection. Its sections are sometimes called \emph{generating sections of symmetries} for the following reasons. Take the restriction $V_{\mathscr E} \to \mathscr E$ of $V$ to $\mathscr E$. Sections of $V_{\mathscr E}$ are locally spanned by the restrictions $v_\alpha|_{\mathscr E}$, $\alpha = 1, \ldots, m$. Take a section $\chi \in \Gamma (V_{\mathscr E})$ and let it be locally given by $\chi = \chi^\alpha v_\alpha|_{\mathscr E}$ for some local functions $\chi^\alpha \in C^\infty (\mathscr E)$. Moreover, let $\tilde \chi \in \Gamma (V)$ be an extension of $\chi$, i.e.{} $\chi = \tilde \chi|_{\mathscr E}$. Finally, let $Y \in \mathfrak X_{\mathscr C}$ be any infinitesimal symmetry of the Cartan distribution on $J^\infty (P, n)$ mapping to $\tilde \chi$ under the projection $X_{\mathscr C} \to \Gamma (V)$ (if we are working locally, we can take, e.g., the evolutionary vector field $\rE_{\tilde \chi}$). If $\mathscr E_0$ is locally given by (\ref{eq:PDE}), one can show that $Y$ is tangent to $\mathscr E$ if and only if locally $\chi$ is a solution of the following \emph{universal linearization} of $\mathscr E_0$:
\[
\ell_{\mathscr E} (\chi) = 0 
\]
 
where 
\[
\ell_{\mathscr E} : \Gamma (V_\mathcal E) \to C^\infty (J^\infty (P,n))^r, \quad r = \operatorname{codim} \mathscr E_0 = \text{number of PDEs}
\]
is the locally defined operator given by
\[
\ell_{\mathscr E}(\chi) = \big( \ell_{\mathscr E}(\chi)_1, \ldots, \ell_{\mathscr E}(\chi)_r \big), \quad \ell_{\mathscr E}(\chi)_a := \left(\rE_{\tilde \chi} F_a\right)|_{\mathscr E} = \frac{\partial F_a}{\partial u^\alpha_I}|_{\mathscr E} D_I|_{\mathscr E} \chi^\alpha.
\]
and    we used that all the total derivatives are tangent to $\mathscr E$  (a more intrinsic and global interpretation of the universal linearization can be provided for PDEs imposed on sections of a fibration, particularly when $\mathcal E_0$ is specified as the zero locus of a \emph{nonlinear differential operator}, see e.g.{} \cite[Chapter 4, Section 2.7]{b...99} or \cite[Section 3.9]{kv98})  . In this case, the restriction $Y|_{\mathscr E}$ is a infinitesimal symmetry of $\mathscr E_0$  

\begin{proposition}
The assignment $\chi \mapsto Y|_{\mathscr E} \operatorname{mod} \Gamma (\mathscr C)$ establishes a bijection between solutions of the universal linearization
\[
\ell_{\mathscr E} (\chi) = 0
\]
of $\mathscr E_0$ and infinitesimal symmetries of $\mathscr E_0$.
\end{proposition}

\begin{remark}
The universal linearization is a linear PDE imposed on sections of the vector bundle $V_{\mathscr E} \to \mathscr E$. As it only involves total derivatives, for any solution $N$ of $\mathscr E_0$, it can be pulled-back along the embedding $j^\infty : N \to \mathscr E$ of Remark \ref{rem:emb_N}. If we do so we get a new linear PDE for sections of the pull-back bundle $j^\infty{}^\ast V$. The latter PDE is just the linearization of $\mathscr E_0$ around the solution $N$. This explain the terminology ``universal linearization''.
\end{remark}

\section{Cohomology of PDEs}\label{Sec:Homol}

\subsection{horizontal cohomology}\label{Sec:hor_cohom}
In this section we attach a cochain complex to every diffiety $(\mathscr E, \mathscr C)$. The associated cohomology is called the \emph{horizontal cohomology} of $\mathscr E$ and it contains important coordinate free information on $\mathscr E_0$.

Let $(\mathscr E, \mathscr C)$ be a diffiety.  By involutivity, the vector subbundle $\mathscr C \to \mathscr E$ of the tangent bundle $T \mathscr E$ is actually a Lie algebroid. Accordingly, there is an associated Cartan calculus. There are also associated De Rham cohomologies (with coefficients). Specifically,    denote by
\[
\overline \Omega{}^\bullet (\mathscr E) := \Gamma (\wedge^\bullet \mathscr C^\ast)
\]
the exterior algebra of the dual of $\mathscr C$. In other words, $\overline \Omega{}^\bullet (\mathscr E)$ consists of differential forms on $\mathscr E$ acting on vector fields in $\mathscr C$. Elements in $\overline \Omega{}^\bullet (\mathscr E)$ are called \emph{horizontal differential forms}. There is a canonical differential $\overline d : \overline \Omega{}^\bullet (\mathscr E) \to \overline \Omega{}^{\bullet +1} (\mathscr E)$, the \emph{horizontal De Rham differential}, defined by the usual Chevalley-Eilenberg formula: for all $\overline \omega \in \overline \Omega{}^q (\mathscr E)$,
\[
\begin{aligned}
& \overline d \overline \omega (X_1, \ldots, X_{q+1}) \\
& = \sum_{i = 1}^{q+1} (-)^{i+1}X_i \left( \overline \omega (X_1, \ldots \widehat{X_i}, \ldots, X_{q+1})\right) \\
& \quad + \sum_{i < j} (-)^{i+j}\overline \omega  \left( [X_i, X_j], \ldots, \widehat{X_i}, \ldots, \widehat{X_j}, \ldots \right), \quad X_1, \ldots, X_{q+1} \in \Gamma (\mathscr C),
\end{aligned}
\]
where a hat denotes omission. The pair $(\overline \Omega{}^\bullet (\mathscr E), \overline d)$ is a commutative differential graded (DG) algebra called the \emph{horizontal De Rham algebra} of $\mathscr E$. Its cohomology is called the \emph{horizontal De Rham cohomology} and denoted $\overline H{}^\bullet (\mathscr E)$. Locally, $\overline \Omega{}^\bullet (\mathscr E)$ is generated, as a graded algebra, by functions and by the \emph{coordinate horizontal $1$-forms}
\[
\overline d x^1, \ldots, \overline d x^n.
\]
In particular, there are no non-trivial horizontal forms of degree higher than $n$, the number of independent variables. In coordinates the horizontal De Rham differential acts as follows:
\[
\overline d \left( f_{i_1 \ldots i_q} \overline d x^{i_1} \wedge \cdots \wedge \overline d x^{i_q} \right) = \left(D_i |_{\mathscr E} f_{i_1 \ldots i_q}\right) \overline d x^i \wedge \overline d x^{i_1} \wedge \cdots \wedge \overline d x^{i_q}.
\]
Notice that the horizontal cohomology is a variant of the standard leaf-wise cohomology of a foliation in our (generically) $\infty$-dimensional setting.

There are various natural DG modules over $(\overline \Omega{}^\bullet (\mathscr E), \overline d)$ coming from standard representations of the Lie algebroid $\mathscr C$. In this subsection we present two of such representations while a third one will be discussed in Section \ref{Sec:Homot}. Denote by
\[
\mathscr V := T\mathscr E / \mathscr C
\]
the \emph{normal bundle} to $\mathscr E$. We will refer to $\mathscr V$ as the \emph{transverse tangent bundle} of $\mathscr E$, and sections of $\mathscr V$ are \emph{transverse vector fileds} (the word ``transverse'' refers to transversality with respect to $\mathscr C$). So, there is a short exact sequence of vector bundles
\begin{equation}\label{eq:SES_C}
0 \longrightarrow \mathscr C \longrightarrow T\mathscr E \longrightarrow \mathscr V \longrightarrow 0.
\end{equation}
The Lie algebroid $\mathscr C$ acts on $\mathscr V$ via the \emph{Bott connection} $\nabla$:
\[
\nabla_X \big(Y \operatorname{mod} \Gamma (\mathscr C) \big) = [X, Y] \operatorname{mod} \Gamma (\mathscr C), \quad \text{for all $X \in \Gamma (\mathscr C)$ and all $Y \in \mathfrak X (\mathscr E)$}.
\]
Accordingly, the graded vector space
\[
\overline \Omega{}^{\bullet} (\mathscr E, \mathscr V) := \Gamma \big(\wedge^\bullet \mathscr C^\ast \otimes \mathscr V\big) \cong \overline \Omega{}^{\bullet} (\mathscr E) \otimes \Gamma (\mathscr V),
\]
is a DG $(\overline \Omega{}^{\bullet} (\mathscr E), \overline d)$-module with differential, denoted 
\[
\overline d : \overline \Omega{}^{\bullet} (\mathscr E, \mathscr V) \to \overline \Omega{}^{\bullet +1 } (\mathscr E, \mathscr V)
\]
as well, given by: for all $\overline W \in \overline \Omega{}^{q} (\mathscr E, \mathscr V)$,
\[
\begin{aligned}
& \overline d\, \overline W (X_1, \ldots, X_{q+1}) \\
& = \sum_{i = 1}^{q+1} (-)^{i+1}\nabla_{X_i} \left( \overline W (X_1, \ldots \widehat{X_i}, \ldots, X_{q+1})\right) \\
& \quad + \sum_{i < j} (-)^{i+j}\overline W  \left( [X_i, X_j], \ldots, \widehat{X_i}, \ldots, \widehat{X_j}, \ldots \right), \quad X_1, \ldots, X_{q+1} \in \Gamma (\mathscr C).
\end{aligned}
\]
The cohomology of the DG module $(\overline \Omega{}^{\bullet} (\mathscr E, \mathscr V), \overline d)$ is denoted $\overline H{}^\bullet (\mathscr E, \mathscr V)$.

Dually to (\ref{eq:SES_C}) there is a short exact sequence of DG algebras
\[
0 \longleftarrow \big(\overline \Omega{}^\bullet (\mathscr E), \overline d\big) \longleftarrow \big( \Omega^\bullet (\mathscr E), d\big) \longleftarrow \big(\mathscr C \Omega^\bullet, d\big) \longleftarrow 0.
\]
Here the projection $\Omega^\bullet (\mathscr E) \to \overline \Omega{}^\bullet (\mathscr E)$ is just the restriction to vector fields in $\mathscr C$, and $\mathscr C \Omega^\bullet$ is its kernel: the $d$-closed graded ideal of forms vanishing when restricted to vector fields in $\mathscr C$. We also denote
\[
\mathscr V \Omega^\bullet = \Gamma (\wedge^\bullet \mathscr V^\ast).
\]
It is a graded subalgebra in $\Omega^\bullet (\mathscr E)$: the graded subalgebra generated by $C^\infty (\mathscr E)$ and $\mathscr C \Omega^1$. We will refer to elements in $\mathscr V \Omega^\bullet$ as \emph{transverse differential forms} on $\mathscr E$. The Lie algebroid $\mathscr C$ acts on $\wedge^\bullet \mathscr V^\ast$ via the \emph{dual Bott connection}, also denoted $\nabla$:
\[
\nabla_X \varpi := \mathcal L_X \varpi, \quad \text{for all $X \in \Gamma (\mathscr C)$ and all $\varpi \in \mathscr V \Omega^\bullet$}.
\]
Accordingly, the graded vector space
\[
\overline \Omega{}^{\bullet} (\mathscr E, \wedge^\bullet \mathscr V^\ast) := \Gamma \big(\wedge^\bullet \mathscr C^\ast \otimes \wedge^\bullet \mathscr V^\ast\big) \cong \overline \Omega{}^{\bullet} (\mathscr E) \otimes \mathscr V^\bullet \Omega,
\]
is also a DG $(\overline \Omega{}^{\bullet} (\mathscr E), \overline d)$-module with differential, denoted 
\[
\overline d : \overline \Omega{}^{\bullet} (\mathscr E, \wedge^\bullet \mathscr V^\ast) \to \overline \Omega{}^{\bullet +1 } (\mathscr E, \wedge^\bullet \mathscr V^\ast),
\]
 given by the usual Chevalley-Eilenberg formula again. The cohomology of the DG module $(\overline \Omega{}^{\bullet} (\mathscr E, \wedge^\bullet \mathscr V^\ast), \overline d)$ is denoted $\overline H{}^\bullet (\mathscr E, \wedge^\bullet \mathscr V^\ast)$.

\subsection{interpreting horizontal cohomologies} Some of the horizontal cohomologies have nice interpretations, sometimes informal, sometimes more precise. In this subsection we illustrate some of them. We begin with horizontal cohomology with trivial coefficients.

\begin{example}[Constrained Variational Principles]\label{exmpl:var_prin}
Top horizontal cohomologies with trivial coefficients $\overline H{}^n (\mathscr E)$ can be (informally) interpreted as \emph{variational principles} constrained by the PDE $\mathscr E_0$. To see this, consider a cohomology class $[L] \in \overline H{}^n (\mathscr E)$. Locally
\[
L = \mathscr L (x, \ldots, u_I, \ldots)\, \overline{\mathrm{d}}{}^n x, \quad \text{where $\overline{\mathrm{d}}{}^n x := \overline d x^1 \wedge \cdots \wedge \overline d x^n$},
\]
for some local function $\mathscr L \in C^\infty (\mathscr E)$. Now, let $N \subseteq P$ be a compact, oriented $n$-dimensional solution of $\mathscr E$ (without boundary). We can pull-back $L$ along the embedding $j^\infty : N \to \mathscr E$. If we do so, we get a honest top form on $N$ that we can integrate using compactness and the orientation. Thus we have a well defined functional
\[
\mathbf F : N \mapsto \mathbf F [N] := \int_N j^\infty{}^\ast L.
\]
If $N$ is locally given by (\ref{eq:N_loc}) and the coordinates $x^1, \ldots, x^n$ on $N$ define the positive orientation, then $\mathbf F$ locally looks like
\begin{equation}\label{eq:F_loc}
\mathbf F [N] = \int \mathscr L \left(x, \ldots, \frac{\partial^{|I|} f}{\partial x^I}, \ldots \right) \mathrm{d}^n x.
\end{equation} 
Actually $\mathbf F$ does only depend on the cohomology class of $L$. Additionally, we see from (\ref{eq:F_loc}) that we can interpret $\mathbf F$, and hence $[L]$ as a (coordinate free version of) a variational principle with Lagrangian density $\mathscr L$, imposed on (certain) solutions of $\mathscr E_0$.
\end{example}

\begin{example}[Conservation Laws]\label{exmpl:cons_law}
Top minus $1$ horizontal cohomologies with trivial coefficients can be interpreted as \emph{conservation laws} for the PDE $\mathscr E_0$. To see this consider a cohomology class $[J] \in \overline H{}^{n-1} (\mathscr E)$. Locally
\[
J = J^i (x, \ldots, u_I, \ldots)\, \overline{\mathrm{d}}{}^{n-1} x_i,
\]
where
\[
\overline{\mathrm{d}}{}^{n-1} x_i := (-)^{i+1}\overline d x^1 \wedge \cdots \wedge \widehat{\overline d x^i} \wedge \cdots \wedge \overline d x^n,
\]
for some local functions $J^i \in C^\infty (\mathscr E)$, $i = 1, \ldots, n$. The cocycle condition $\overline d J = 0$ in local coordinates reads
\[
D_i |_{\mathscr E} J^i = 0,
\]
so $J$ can be seen as a \emph{conserved current} along solutions of $\mathscr E_0$. Now, let $N \subseteq P$ be a solution of $\mathscr E$ and let $\Sigma \subseteq N$ be a compact, oriented hypersurface (without boundary) in $N$. We can pull-back $J$ along $j^\infty : N \to \mathscr E$. If we do so, we get a honest $(n-1)$-form on $N$ that we can integrate on $\Sigma$ using compactness and the orientation. The integral
\[
\int_\Sigma j^\infty{}^\ast J
\]
does only depend on (the horizontal cohomology class $[J]$ and on) the homology class of $\Sigma$ in $N$. In this sense, it is conserved along every $N$. Summarizing, $[J]$ can be seen as the conservation law associated to the conserved current $J$. Notice that lower degree horizontal cohomologies can be similarly seen as lower degree conserved charges to be integrated on higher codimensional submanifolds. For instance the electric charge is a degree $(n-2)$ horizontal cohomology class for the Maxwell equations on a vacuous $n$-dimensional space-time.
\end{example}

\begin{example}[Infinitesimal Symmetries]\label{exmpl:inf_sym}
In this example we discuss degree $0$ horizontal cohomologies with coefficients in $\mathscr V$. The $0$-th cohomology $\overline H{}^0 (\mathscr E, \mathscr V)$ is just the kernel of the operator
\[
\overline d : \Gamma (\mathscr V) \to \overline \Omega{}^1 (\mathscr E, \mathscr V)
\]
and it consists of \emph{parallel sections} with respect to the Bott connection. It is easy to see that a section $Y \operatorname{mod} \Gamma (\mathscr C)$ of $\mathscr V$, with $Y \in \mathfrak X (\mathscr E)$, is parallel with respect to the Bott connection if and only if $Y$ is an infinitesimal symmetry of the $\mathscr E_0$. We conclude that the natural inclusion
\[
\mathfrak X_\mathscr{E} / \Gamma (\mathscr C) \to \mathfrak X (\mathscr E) / \Gamma (\mathscr C) = \Gamma (\mathscr V)
\]
establishes a bijection between non-trivial infinitesimal symmetries of $\mathscr E_0$ and the $0$-th cohomolgy $\overline H{}^0 (\mathscr E, \mathscr V)$. 

Let $Y \in \mathfrak X_{\mathscr C}$.  Remember from Remark \ref{rem:infty_dim} that the vector field $Y$ might not possess a flow. However,    when $Y$ possesses a flow, then the latter maps $n$-dimensional integral submanifolds of $\mathscr C$ to $n$-dimensional integral submanifolds. So, morally, it defines a flow on the space of solutions of $\mathscr E_0$. If two infinitesimal symmetries differ by a trivial infinitesimal symmetry then their flows (if they exist) induce the same flow on the space of solutions, because the flow of a trivial symmetry fixes every $n$-dimensional integral submanifold. We conclude that nontrivial infinitesimal symmetries can be informally seen as vector fields on the space of solutions. This remark will be relevant in the next subsection.
\end{example}

\begin{example}[Euler Lagrange Equations]\label{exmpl:ELeq}
We now discuss top horizontal cohomology with coefficients in transverse $1$-forms $\mathscr V^\ast$. We begin noticing that there is a map
\[
d^{\mathrm{sec}} : \overline H{}^n (\mathscr E) \to \overline H{}^n (\mathscr E, \mathscr V^\ast)
\]
defined as follows. Let $[L] \in \overline H{}^n (\mathscr E)$, and choose any $n$-form $\omega \in \Omega^n (\mathscr E)$ projecting to $L \in \overline \Omega{}^\bullet (\mathscr E)$ under $\Omega^\bullet (\mathscr E) \to \overline \Omega{}^\bullet (\mathscr E)$. The differential $d \omega$ belongs to the ideal $\mathscr C \Omega^\bullet$. Hence, it has the following property: the contraction
\[
i_{X_n} \cdots i_{X_1} d\omega 
\]
with any $n$ vector fields $X_1, \ldots, X_n$ in $\mathscr C$ is a transverse differential $1$-form: 
\[
i_{X_n} \cdots i_{X_1} d\omega \in \mathscr V \Omega^1 = \mathscr C \Omega^1.
\]
 Hence we have a well-defined, skew-symmetric and $C^\infty (\mathscr E)$-multilinear map
\[
\varpi : \Gamma (\mathscr C) \times \cdots \times \Gamma (\mathscr C) \to \mathscr V \Omega^1, \quad (X_1, \ldots, X_n) \mapsto i_{X_n} \cdots i_{X_1} d\omega 
\]
which we can interpret as an element in $\overline \Omega{}^n (\mathscr C, \mathscr V^\ast)$, also denoted  $\varpi$. As $dd\omega = 0$, we get $\overline d \varpi = 0$. The cohomology class $[\varpi ] \in \overline H{}^n (\mathscr E, \mathscr V^\ast)$ does only depend on the cohomology class $[L]$, and we put
\[
d^{\mathrm{sec}} [L] = [ \varpi]. 
\]
We now provide a more concrete description of the map $d^{\mathrm{sec}} : \overline H{}^n (\mathscr E) \to \overline H{}^n (\mathscr E, \mathscr V^\ast)$. For simplicity, from now on, in this example, we will only consider the case when $\mathscr E_0$ is the empty PDE $0 = 0$, hence $\mathscr E = J^\infty (P, n)$. For the sake of brevity, we denote $J^\infty = J^\infty (P, n)$. When $\mathscr E = J^\infty$, there is a simple description of $\overline H{}^n (\mathscr E, \mathscr V^\ast)$. Namely, there is a vector space isomorphism
\[
\Gamma \big(V^\ast \otimes \wedge^n \mathscr C^\ast \big) \cong \overline H{}^n (J^\infty, \mathscr V^\ast)
\]
defined as follows. First notice that there is a vector bundle projection
\[
\mathscr V \to V
\]
mapping a transverse vector $v \operatorname{mod} \mathscr C$ at the point $z = N_p^\infty \in J^\infty$ to 
\[
d\pi (v) \operatorname{mod} T_p N \in T_p P / T_p N = V_z.
\]
Dually, there is a vector bundle inclusion
\[
V^\ast \hookrightarrow \mathscr V^\ast.
\]
Accordingly, there is an inclusion
\begin{equation}\label{eq:incl_V*}
\Gamma \big(V^\ast \otimes \wedge^n \mathscr C^\ast \big) \hookrightarrow \overline \Omega{}^n (J^\infty, \mathscr V^\ast).
\end{equation}
that can be described in local coordinates as follows. Let $v^1, \ldots, v^m \in \Gamma (V^\ast)$ be the dual basis of the local basis $v_1, \ldots, v_m \in \Gamma (V)$ described in Subsection \ref{Sec:CD}. Then the inclusion (\ref{eq:incl_V*}) maps
\[
v^\alpha \otimes \overline{\mathrm d}{}^n x \mapsto \overline{\mathrm d}{}^n x \otimes (du^\alpha - u^\alpha_i dx^i), \quad \alpha = 1, \ldots, m.
\]

\begin{proposition}
The composition
\begin{equation}\label{eq:isom_top_V*}
\Gamma \big(V^\ast \otimes \wedge^n \mathscr C^\ast \big) \hookrightarrow \overline \Omega{}^n (J^\infty, \mathscr V^\ast) \longrightarrow \overline H{}^n (J^\infty, \mathscr V^\ast)
\end{equation}
is a vector space isomorphism. Moreover, for $L \in \overline \Omega{}^n (J^\infty)$ locally given by 
\[
L = \mathscr L (x, \ldots, u_I, \ldots)\, \overline{\mathrm d}{}^n x
\]
the section $\mathbf E [L]$ of $V^\ast \otimes \wedge^n \mathscr C^\ast$ corresponding to $d^{\mathrm{sec}} [L]$ via (\ref{eq:isom_top_V*}) is locally given by
\[
\mathbf E [L] = \sum_{|I| < \infty} (-)^{|I|} D_I \left( \frac{\partial \mathscr L}{\partial u^\alpha_I}\right) v^\alpha \otimes \overline{\mathrm d}{}^n x.
\]
\end{proposition}

In view of the first part of the above proposition, the top cohomology $\overline H{}^n (J^\infty, \mathscr V^\ast)$ is in bijection with the $C^\infty (J^\infty)$-module of sections of a vector bundle over $J^\infty$. Hence it makes sense to speak about the \emph{zero locus} in $J^\infty$ of a cohomology class in $\overline H{}^n (J^\infty, \mathscr V^\ast)$. From the second part of the proposition, the zero locus of $d^{\mathrm{sec}} [L]$, equivalently the zero locus of $\mathbf E [L]$, is the \emph{Euler-Lagrange equations} corresponding to the variational principle $[L]$ (see Example \ref{exmpl:var_prin}):
\[
\sum_{|I| < \infty} (-)^{|I|} D_I \left( \frac{\partial \mathscr L}{\partial u^\alpha_I}\right)  = 0. 
\]
Similarly, one can show that there is a map
\[
d^{\mathrm{sec}} : \overline H{}^n (J^\infty, \mathscr V^\ast) \to \overline H{}^n (J^\infty, \wedge^2 \mathscr V^\ast)
\]
(see the next subsection) mapping a system of $m$ PDEs
\[
 E_\alpha (x, \ldots, u_I, \ldots) = 0 
\]
to the associated \emph{Helmoltz conditions} for variationality.
\end{example}

\subsection{algebraic structures on horizontal cohomologies}\label{Subsec:LR}
 Examples \ref{exmpl:var_prin} and \ref{exmpl:cons_law} suggest that the horizontal cohomology $\overline H{}^\bullet (\mathscr E)$ should be interpreted as ``smooth functions on the space of solutions of the PDE $\mathscr E_0$''. For this reason we also denote it
 \[
 \boldsymbol C^\infty := \overline H{}^\bullet (\mathscr E)
 \] 
 and call it \emph{secondary functions} adopting a terminology by Vinogradov. Similarly, the discussion at the end of Example \ref{exmpl:inf_sym} suggests that $\overline H{}^\bullet (\mathscr E, \mathscr V)$ should be interpreted as ``vector fields on the space of solutions of $\mathscr E_0$''. For this reason we also denote it
 \[
 \boldsymbol{\mathfrak X} := \overline H{}^\bullet (\mathscr E, \mathscr V)
 \]
 and call it \emph{secondary vector fields}. Finally, Example \ref{exmpl:ELeq} suggests that $\overline H{}^\bullet (\mathscr E, \mathscr V)$ should be interpreted as ``differential forms on the space of solutions of $\mathscr E_0$'', we denote
 \[
 \boldsymbol{\Omega}^p := \overline H{}^\bullet (\mathscr E, \wedge^p \mathscr V^\ast)
 \]
 and call it \emph{secondary differential $p$-forms}, $p = 0, 1, \ldots$. This interpretation is supported by the fact that the above cohomologies possess the appropriate algebraic structures for functions, vector fields, differential forms. For instance, $\boldsymbol C^\infty$ is naturally a (graded) commutative algebra, being the cohomology of a  commutative DG algebra. Similarly, the pair $(\boldsymbol C^\infty, \boldsymbol{\mathfrak X})$ is a (graded) Lie-Rinehart algebra. We now explain the latter claim.
 
 Choose once for all a splitting of the short exact sequence (\ref{eq:SES_C}):
 \begin{equation}\label{eq:split}
 \xymatrix{0 \ar[r] & \mathscr C \ar[r] & T\mathscr E \ar[r] & \mathscr V \ar[r] \ar@/^9pt/[l]& 0
 }.
 \end{equation}
 
\vspace{10pt}
\noindent This induces a direct sum decomposition $T\mathscr E \cong \mathscr C \oplus \mathscr V$ that allows us to interpret $\mathscr V$ as a distribution and transverse vector fields as honest vector fields on $\mathscr E$. Dually, we get a factorization
\begin{equation}\label{eq:fact}
\Omega^\bullet (\mathscr E) \cong \overline \Omega{}^{\bullet} (\mathscr E) \otimes \mathscr V \Omega^\bullet
\end{equation}
that allows us to interpret horizontal differential forms as honest differential forms on $\mathscr E$, and in the rest of this section we will do so.  Similarly, we understand the embedding $\overline \Omega{}^\bullet (\mathscr E, \mathscr V) \hookrightarrow \Omega^\bullet (\mathscr E, T\mathscr E)$ induced by the splitting (\ref{eq:split}) (and the factorization (\ref{eq:fact})) and interpret $\mathscr V$-valued horizontal differential forms as honest vector valued forms on $\mathscr E$.    Now, recall that vector valued forms are naturally equipped with a graded Lie bracket: the Fr\"olicher-Nijenhuis bracket $[-,-]^{\mathrm{fn}}$ (see, e.g., \cite[Section 16]{m08}). Moreover, the resulting graded Lie algebra acts on differential forms via the Lie derivative (along vector valued forms). 

Notice that, unless the distribution $\mathscr V \subseteq T \mathscr E$ is itself involutive, the Fr\"olicher-Nijenhuis bracket does not preserve $\mathscr V$-valued horizontal forms $\overline \Omega{}^\bullet (\mathscr E, \mathscr V) $. We can define a bracket on $\overline \Omega{}^\bullet (\mathscr E, \mathscr V)$ by first taking the Fr\"olicher-Nijenhuis bracket, and then projecting to $\overline \Omega{}^\bullet (\mathscr E, \mathscr V)$, but the bracket obtained in this way is not a Lie bracket (unless $\mathscr V$ is involutive). However, it gives a honest graded Lie bracket on $\boldsymbol{\mathfrak X}$. Namely, denote by
\[
\operatorname{pr} : \Omega^\bullet (\mathscr E) \to \overline \Omega{}^\bullet (\mathscr E), \quad\quad \operatorname{pr} : \Omega^\bullet (\mathscr E , T \mathscr E) \to \overline \Omega{}^\bullet (\mathscr E, \mathscr V)
\]
the (natural) projections. Define a bracket $[-,-]^{\mathrm{sec}}$ on $\boldsymbol{\mathfrak X}$, the \emph{secondary commutator}, via:

\begin{equation}\label{eq:sec_Lie_bra}
[-,-]^{\mathrm{sec}} : \boldsymbol{\mathfrak X} \times \boldsymbol{\mathfrak X} \to \boldsymbol{\mathfrak X}, \quad \left(\boldsymbol W, \boldsymbol U \right) \mapsto \left[ \boldsymbol W, \boldsymbol U \right]^{\mathrm{sec}} := \big[\operatorname{pr}[\overline W, \overline U]^{\mathrm{fn}}\big],
\end{equation}
where $\boldsymbol W = [\overline W], \boldsymbol U = [\overline U]$ are the cohomology classes of the cocycles $\overline W, \overline U \in \overline \Omega{}^\bullet (\mathscr E, \mathscr V)$ respectively.

We also define the \emph{secondary Lie derivative}:
\begin{equation}\label{eq:sec_Lie_der}
\mathcal L^{\mathrm{sec}} : \boldsymbol{\mathfrak X} \times \boldsymbol{C}^\infty \to  \boldsymbol{C}^\infty, \quad \left(\boldsymbol W, \boldsymbol f \right) \mapsto \mathcal L^{\mathrm{sec}}_{\boldsymbol W} \boldsymbol f :=\big[ \operatorname{pr} \mathcal L_{\overline W} \overline \omega \big],
\end{equation}
where $\boldsymbol f = [\overline \omega]$ is the cohomology class of the cocycle $\overline \omega \in \overline \Omega{}^\bullet (\mathscr E)$. 

\begin{theorem}\label{theor:LR}
The operations (\ref{eq:sec_Lie_bra}) and (\ref{eq:sec_Lie_der}) are well-defined and independent of the splitting (\ref{eq:split}). Together with them and the graded $\boldsymbol C^\infty$-module structure on $\boldsymbol{\mathfrak X}$ induced by the DG $(\overline \Omega{}^\bullet (\mathscr E), \overline d)$-module structure on $(\overline \Omega{}^\bullet (\mathscr E, \mathscr V), \overline d)$, the pair $(\boldsymbol{C}^\infty, \boldsymbol{\mathfrak X})$ form a graded Lie-Rinehart algebra.
\end{theorem}

As for secondary differential forms, they form a commutative DG algebra. To see this, consider the splitting (\ref{eq:split}) again and understand the induced factorization (\ref{eq:fact}). For every $p, p'$, we have a \emph{secondary wedge product}:
\begin{equation}\label{eq:wedge_sec}
\wedge^{\mathrm{sec}} : \boldsymbol \Omega^p \times \boldsymbol \Omega^{p'} \to \boldsymbol \Omega^{p+p'}, \quad (\boldsymbol \omega, \boldsymbol \omega') \mapsto \boldsymbol \omega \wedge^{\mathrm{sec}}  \boldsymbol \omega' := \big[ \omega \wedge \omega' \big],
\end{equation}
where $\boldsymbol \omega = [\omega], \boldsymbol \omega' = [\omega']$ are the cohomology classes of the cocycles $\omega \in \overline \Omega{}^\bullet (\mathscr E, \wedge^p \mathscr V^\ast)$, $\omega' \in \overline \Omega{}^\bullet (\mathscr E, \wedge^{p'} \mathscr V^\ast)$.

We also have a \emph{secondary differential}.  In order to define it we first notice that the De Rham differential does not map $\wedge^p \mathscr V^\ast$-valued horizontal form to $\wedge^{p+1} \mathscr V^\ast$-valued horizontal form (unless $\mathscr V$ is involutive). We can define a map
\[
\overline \Omega{}^\bullet (\mathscr E, \wedge^p \mathscr V^\ast ) \to \overline \Omega{}^\bullet (\mathscr E, \wedge^{p+1} \mathscr V^\ast )
\]
by first taking the De Rham differential, and then projecting to $\Omega{}^\bullet (\mathscr E, \wedge^{p+1} \mathscr V^\ast )$, but the maps obtained in this way do not square to zero (unless $\mathscr V$ is involutive). However, they give a honest differential on $\boldsymbol{\Omega}^\bullet$. Namely, denote by
\[
\operatorname{pr} : \Omega^\bullet (\mathscr E) \to \Omega{}^\bullet (\mathscr E, \wedge^{p+1} \mathscr V^\ast)
\]
the projection. Define a map $d^{\mathrm{sec}}$ on $\boldsymbol{\Omega}^\bullet$ via:
\begin{equation}\label{eq:d_sec}
d^{\mathrm{sec}} : \boldsymbol \Omega^p \to \boldsymbol \Omega^{p+1}, \quad \boldsymbol \omega \mapsto d^{\mathrm{sec}} \boldsymbol \omega := \big[\operatorname{pr} d \omega \big].
\end{equation}

\begin{theorem}\label{theor:DGA}
The operations (\ref{eq:wedge_sec}) and (\ref{eq:d_sec}) are well-defined and independent of the splitting (\ref{eq:split}). Together with them $\boldsymbol{\Omega}^\bullet$ form a commutative DG algebra (with respect to the \emph{total grading}, i.e.{} the grading on $\overline H{}^q (\mathscr E, \wedge^p \mathscr V^\ast)$ given by $p+q$).
\end{theorem}

\begin{remark}
The secondary differential (\ref{eq:d_sec}) agrees with that introduced in Example \ref{exmpl:ELeq}.
\end{remark}

There is a more conceptual way to introduce the operations (\ref{eq:wedge_sec}) and (\ref{eq:d_sec}), based on the \emph{Vinogradov's $\mathscr C$-spectral sequence}, that we now explain. Namely, recall the $d$-closed graded ideal $\mathscr C \Omega^\bullet \subseteq \Omega^\bullet (\mathscr E)$ from Subsection \ref{Sec:hor_cohom} (it consists of differential forms vanishing when restricted to vector fields in $\mathscr C$).  The powers 
\[
\mathscr C^p \Omega^\bullet := \underset{\text{$p$ times}}{\underbrace{\mathscr C \Omega^\bullet \wedge \cdots \wedge \mathscr C \Omega^\bullet}},\quad p = 0, 1, \ldots
\]
define a filtration
\begin{equation}\label{eq:filtr}
\Omega^\bullet (\mathscr E) = \mathscr C^0 \Omega^\bullet \supseteq \mathscr C \Omega^\bullet \supseteq \cdots \supseteq \mathscr C^p \Omega^\bullet \supseteq \cdots, 
\end{equation}
by graded $d$-closed ideals, which in turn determines a spectral sequence $\mathscr C E$, called the \emph{$\mathscr{C}$-spectral sequence}, computing the ordinary De Rham cohomologies of $\mathscr E$.

\begin{proposition}
For all $p, q$, the map
\[
\overline \Omega{}^\bullet (\mathscr E, \wedge^p \mathscr V^\ast) \to \mathscr E_0^{p,\bullet} = \mathscr C^p \Omega^{p+\bullet} / \mathscr C^{p+1} \Omega^{p+\bullet}, \quad \overline \omega \otimes \varpi \mapsto \omega \wedge \varpi \operatorname{mod} \mathscr C^{p+1} \Omega^{p+\bullet}
\]
is a well-defined cochain isomorphism, where $\overline \omega \in \overline \Omega{}^\bullet (\mathscr E)$, $\varpi \in \mathscr V \Omega^p$, and $\omega \in \Omega^\bullet (\mathscr E)$ is any $q$-form projecting to $\overline \omega$ under $\Omega^\bullet (\mathscr E) \to \overline \Omega{}^\bullet (\mathscr E)$. The induced isomorphism
\[
\boldsymbol \Omega^p = \overline H{}^\bullet (\mathscr E, \wedge^p \mathscr V^\ast) \to \mathscr C E^{p, \bullet}_1 = H^\bullet (\mathscr C E^{p, \bullet}_0, d_0^{p, \bullet})
\]
identifies $d^{sec}$ with the differential $d_1$ in the first page of the $\mathscr C$-spectral sequence, and $\wedge^{sec}$ with the product structure on $\mathscr C E_1$ given by the fact that (\ref{eq:filtr}) is a filtration by DG ideals.
\end{proposition}

Hence $\wedge^{\mathrm{sec}}, d^{\mathrm{sec}}$ could have been defined via the $\mathscr C$-spectral sequence without using the splitting (\ref{eq:split}). One can even show that there is an analogue of the whole Cartan Calculus on secondary differential forms and secondary vector fields. All these operations together form the fundamentals of \emph{Vinogradov's Secondary Calculus}. Most of Vinogradov's work on PDEs has been based on the following principle: \emph{Secondary Calculus is an appropriate replacement for Differential Calculus on the space of solutions of a PDE}.

We conclude this sections presenting an informal \emph{secondarization recipe}. Namely, let $\Phi$ be any natural construction in differential geometry (e.g.{} functions on manifolds, vector fields, the De Rham complex, etc.). The discussion in this subsection suggests the following recipe to find the \emph{secondary version of $\Phi$}, i.e.{} an appropriate replacement $\boldsymbol \Phi$ for the construction $\Phi$ on the space of solutions of a PDE $\mathscr E_0$.

\begin{recipe} \quad
\begin{enumerate}
\item Define the \emph{transverse version} $\mathscr V \Phi$ of $\Phi$ (when $\Phi$ is ``vector fields'', then its transverse version is $\Gamma (\mathscr V)$).
\item Notice that $\mathscr V \Phi$ is canonically equipped with an action of the Lie algebroid $\mathscr C$ giving a DG $(\overline \Omega{}^\bullet(\mathscr E) , \overline d)$-module structure to $\overline \Omega{}^\bullet(\mathscr E) \otimes \mathscr V \Phi$ (when $\Phi$ is ``vector fields'', then $\overline \Omega{}^\bullet(\mathscr E) \otimes \mathscr V \Phi$ is $\overline \Omega{}^\bullet (\mathscr E, \mathscr V)$).
\item Define $\boldsymbol \Phi$ as the horizontal cohomology with coefficients in $\mathscr V \Phi$: $\boldsymbol \Phi := H^\bullet (\overline \Omega(\mathscr E) \otimes \mathscr V \Phi, \overline d)$ (when $\Phi$ is ``vector fields'', then $\boldsymbol \Phi$ is $\boldsymbol{\mathfrak X}$).
\item Find the appropriate algebraic structures on $\boldsymbol \Phi$ (when $\Phi$ is ``vector fields'', then $\boldsymbol \Phi = \boldsymbol{\mathfrak X}$ is a graded Lie-Rinehart algebra over $\boldsymbol C^\infty$).
\end{enumerate}
\end{recipe}

In the next section we will use this recipe to define \emph{secondary scalar differential operators}, i.e.{} an appropriate replacement for scalar differential operators on the space of solutions of a PDE.

\section{Homotopy of PDEs}\label{Sec:Homot}

In this section we quickly illustrate the latest developments on Vinogradov's ideas about Secondary Calculus. They are not due to Vinogradov himself but they are definitely inspired by his work.

\subsection{homotopy algebras}

In this subsection, we recall the definitions of those homotopy algebras that we will need in the rest of the paper. Roughly, a homotopy algebra of a certain type (associative, Lie, etc.) is a graded vector space equipped with operations satisfying the algebra axioms only up to a coherent system of \emph{higher homotopies}. Homotopy algebras appear when cohomology and homotopy interact with algebraic structures \cite{lv12, v14}. We adopt the degree $1$ convention on the operations of a homotopy algebra, and we only work with (graded) vector spaces over a field of zero characteristic.

\begin{definition}
An \emph{$L_\infty$-algebra} \cite{ls93, lm95} is a graded vector space $V$ equipped with a sequence $\mathfrak l = \{ \mathfrak l_k \}_{k \in \mathbb N}$ of degree $1$ symmetric multilinear brackets:
\[
\mathfrak l_k : \left(S^k V\right){}^\bullet \to V^{\bullet +1}, \quad (v_1, \ldots, v_{k}) \mapsto \mathfrak l_k (v_1, \ldots, v_{k})
\]
satisfying the following \emph{higher Jacobi identities}:
\begin{equation}\label{eq:high_Jac}
\sum_{r+s = k} \sum_{\sigma \in S_{r, s}}\epsilon (\sigma , v) \mathfrak l_{s+1} \big(\mathfrak l_r (v_{\sigma(1)}, \ldots, v_{\sigma(r)}), v_{\sigma (r+1)}, \ldots, v_{\sigma (r+s)}\big) = 0,
\end{equation}
for all $k \in \mathbb N$, where $S_{r,s}$ denotes $(r,s)$-unshuffles, and $\epsilon (\sigma, v)$ is the Koszul sign associated to the permutation $\sigma$ of the homogeneous elements $v_1, \ldots, v_k \in V$.
\end{definition}

For $k = 1$ the identity (\ref{eq:high_Jac}) says that $\mathfrak l_1 : V^\bullet \to V^{\bullet +1}$ is a differential. In particular there is a cohomology $H^\bullet (V, \mathfrak l_1)$ for every $L_\infty$-algebra $(V , \mathfrak l)$. For $k = 2$ the identity (\ref{eq:high_Jac}) says that $\mathfrak l_1$ is a derivation with respect to the binary bracket $\mathfrak l_2 : (S^2 V)^\bullet \to V^{\bullet +1}$. For $k = 3$ the identity (\ref{eq:high_Jac}) says that $\mathfrak l_2$ is a graded Lie bracket only up to a homotopy encoded by $\mathfrak l_3$. In particular $\mathfrak l_2$ induces a honest graded Lie bracket in cohomology $H^\bullet (V, \mathfrak l_1)$ (up to a d\'ecalage isomorphism that we will ignore, for simplicity). So, whenever the cohomology of a cochain complex supports a graded Lie bracket, it is natural to wonder whether or not it comes from an $L_\infty$-algebra structure on cochains.

There is also a notion of $L_\infty$-module over an $L_\infty$-algebra.

\begin{definition}
An \emph{$L_\infty$-module} \cite{lm95} over an $L_\infty$-algebra $(V, \mathfrak l)$ is a graded vector space $W$ equipped with a sequence $\mathfrak m = \{ \mathfrak m_k \}_{k \in \mathbb N}$ of degree $1$ multilinear brackets:
\[
\mathfrak m_k : \left(S^{k-1} V \otimes W \right){}^\bullet \to W^{\bullet + 1}, \quad (v_1, \ldots, v_{k-1}, w) \mapsto \mathfrak m_k (v_1, \ldots, v_{k-1}| w)
\]
satisfying:
\[
\begin{aligned}
& \sum_{r+s = k-1} \sum_{\sigma \in S_{r, s}}\epsilon (\sigma , v) \Big[\mathfrak m_{s+1} \big(\mathfrak l_r (v_{\sigma(1)}, \ldots, v_{\sigma(r)}), v_{\sigma (r+1)}, \ldots, v_{\sigma (r+s)}| w \big)  \\
& + (-)^{|v_{\sigma(1)}| + \cdots +|v_{\sigma(r)}|} \mathfrak m_{r+1} \big( v_{\sigma (1)}, \ldots, v_{\sigma (r)}| \mathfrak m_s (v_{\sigma(r+1)}, \ldots, v_{\sigma(r+s)}|w)\big) \Big]= 0
\end{aligned}
\]
for all $k \in \mathbb N$, $v_1, \ldots, v_{k-1} \in V$, $w \in W$.
\end{definition}

In particular $\mathfrak m_1 : W^\bullet \to W^{\bullet +1}$ is a differential and $\mathfrak m_2$ induces a honest graded $H^\bullet (V, \mathfrak l_1)$-module structure on $H^\bullet (W, \mathfrak m_1)$.

The next definition is a homotopy version of the notion of Lie-Rinehart algebra (see \cite{k01, k01b, h05, h13, vi12, v15}).

\begin{definition}
An \emph{$LR_\infty$-algebra} \cite{vi12} (\emph{LR} for Lie-Rinehart) is a pair $(A, L)$ consisting of two graded vector spaces equipped with the following additional structures:
\begin{enumerate}
\item $L$ is an $L_\infty$-algebra with structure maps $\mathfrak l$;
\item $A$ is a commutative DG algebra with differential denoted $\mathfrak m_1$;
\item $(L, \mathfrak l_1)$ is a DG $(A^\bullet, \mathfrak m_1)$-module;
\item $A$ is an $L_\infty$-module over $(L^\bullet, \mathfrak l)$ with structure maps $\mathfrak m$ (this means that the $1$-ary bracket is precisely $\mathfrak m_1$).
\end{enumerate} 
Additionally, the maps
\[
\mathfrak m_k : \left(S^{k-1} L \otimes A \right){}^\bullet \to A^{\bullet + 1}
\]
are graded $A$-multilinear in the first $(k-1)$-arguments and they are a derivation in the last argument. Finally the $L_\infty$-algebra brackets $\mathfrak l$ satisfy the following \emph{Leibniz rule}:
\[
\begin{aligned}
&\mathfrak l_k (v_1, \ldots, v_{k-1}, a v_k) \\
& = \mathfrak m_{k} (v_1, \ldots, v_{k-1}|a) v_k + (-)^{|a| \left(|v_1| + \cdots + |v_{k-1}| + 1 \right)}a \mathfrak l_k (v_1, \ldots, v_{k-1}, v_k),
\end{aligned}
\]
for all $v_1, \ldots, v_k \in L$, and $a \in A$.
\end{definition}

\begin{remark}\label{rem:LR_cohom}
Let $(A, L)$ be an $LR_\infty$-algebra with structure maps $(\mathfrak m, \mathfrak l)$. Similarly as for $L_\infty$-algebras and $L_\infty$-modules, the cohomology $(H^\bullet (A, \mathfrak m_1), H^\bullet (L, \mathfrak l_1))$ is then a graded Lie-Rinehart algebra.
\end{remark}

As for Lie algebroids, there is a \emph{Chevalley-Eilenberg construction} for $LR_\infty$-algebras. Namely, let $(A, L)$ be an $LR_\infty$-algebra. Consider the graded commutative algebra $\operatorname{Sym}^\bullet(A, L)$ consisting of graded symmetric $A$-multilinear maps
\[
S^\bullet L \to A.
\]
The structure maps $(\mathfrak l, \mathfrak m)$ of $(A, L)$ induce on $\operatorname{Sym}^\bullet (A, L)$ a sequence $\mathfrak d = \{ \mathfrak d_k\}$ of degree $1$ derivations:
\begin{equation}\label{eq:d_Sym}
\mathfrak d_k : \operatorname{Sym}^h (A, L) \to \operatorname{Sym}^{k+h-1} (A, L)
\end{equation}
given by
\[
\begin{aligned}
&\mathfrak d_k \omega (v_1, \ldots, v_{h}) \\
& = \sum_{\sigma \in S_{k-1,h}} \epsilon(\sigma, v) (-)^{|\omega| \sum_{i = 1}^{k-1}|v_{\sigma(i)}|}
\mathfrak m_k \big(v_{\sigma(1)}, \ldots, v_{\sigma(k-1)}|\omega (v_{\sigma (k)}, \ldots, v_{\sigma (k+h-1)})\big) \\
& \quad + \sum_{\sigma \in S_{k, h-1}} \epsilon(\sigma, v) \omega \big(\mathfrak l_k (v_{\sigma(1)}, \ldots, v_{\sigma (k)}), v_{\sigma (k+1)}, \ldots, v_{\sigma (k+h-1)}\big).
\end{aligned}
\]

\begin{proposition}
The maps $\mathfrak d_k$ are well-defined degree $1$ derivations of the graded algebra $\operatorname{Sym}^\bullet(A, L)$. Additionally they satisfy
\begin{equation}\label{eq:d_d}
\sum_{r+s = k} \big[\mathfrak d_r, \mathfrak d_s \big] = 0
\end{equation}
for all $k \in \mathbb N$. In particular the total derivation
\[
\mathfrak D := \mathfrak d_1 + \mathfrak d_2 + \cdots
\]
squares to zero, whenever defined. 
\end{proposition}

The pair $(\operatorname{Sym}^\bullet(A, L), \mathfrak d)$ is called the \emph{Chevalley-Eilenberg} algebra of $(A, L)$.

\begin{remark}\label{rem:DGA_cohom}
The derivation $\mathfrak d_2$ does not square to $0$ in general (it does only up to a homotopy encoded by $\mathfrak d_3$), however it induces a derivation which squares to zero in the cohomology of $\mathfrak d_1$. Hence 
$
H^\bullet \left( \operatorname{Sym}(A, L), \frak d_1\right)
$, with the derivation induced by $\mathfrak d_2$, is a commutative DG algebra.
\end{remark}

\begin{remark}
Let $(A, L)$ be an $LR_\infty$-algebra. If $L$ is projective and finitely generated as an $A$-module, then the Chevalley-Eilenberg algebra knows everything about $(A, L)$. More precisely, let $A$ be a commutative graded algebra and let $L$ be a projective and finitely generated graded $A$-module. Then the Chevalley-Eilenberg construction establishes a bijection between $LR_\infty$-algebra structures on $(A, L)$ extending the preexisting structures on one side, and sequences $\mathfrak d = \{ \mathfrak d_k\}$ of degree $1$ derivations of $\operatorname{Sym}^\bullet(A, L)$ satisfying both (\ref{eq:d_Sym}) and (\ref{eq:d_d}) on the other side.
\end{remark}

\subsection{the $LR_\infty$-algebra of secondary vector fields}

Let $(\mathscr E, \mathscr C)$ be a diffiety. According to Theorem \ref{theor:LR} the pair $(\boldsymbol C^\infty, \boldsymbol{\mathfrak X})$ is a graded Lie-Rinehart algebra. Recall that $\boldsymbol C^\infty = \overline H{}^\bullet (\mathscr E) = H^\bullet (\overline \Omega (\mathscr E), \overline d)$ and $\boldsymbol{\mathfrak X} = \overline H{}^\bullet (\mathscr E, \mathscr V) = H^\bullet(\overline \Omega(\mathscr E, \mathscr V), \overline d)$. Moreover, the cochains $\overline \Omega{}^\bullet(\mathscr E, \mathscr V)$ are naturally a graded module over the cochains $\overline \Omega{}^\bullet (\mathscr E)$. It is then natural to look for an $LR_\infty$-algebra structure on cochains $(\overline \Omega{}^\bullet (\mathscr E), \overline \Omega{}^\bullet(\mathscr E, \mathscr V))$ responsible for the Lie-Rinehart algebra structure on $(\boldsymbol C^\infty, \boldsymbol{\mathfrak X})$. Such $LR_\infty$-algebra structure exists indeed as we now show. 

Choose again a splitting (\ref{eq:split}) and use it to see transverse vector fields and horizontal forms as honest vector fields and differential forms on $\mathscr E$, as we did in Subsection \ref{Subsec:LR}. We will also interpret horizontal forms with values in transverse forms as ordinary vector valued forms.

The normal bundle $\mathscr V$ is now a (non-necessarily involutive) distribution on $\mathscr E$. Its curvature is 
\[
R : \wedge^2 \mathscr V \to \mathscr C, \quad (Y, Z) \mapsto R(Y, Z) := [Y, Z] - \operatorname{pr} [Y, Z],
\]  
where we denote by $\operatorname{pr} : T\mathscr E \to \mathscr V$ the projection, and can be seen as a vector valued $2$-form on $\mathscr E$:
\[
R \in \mathscr V\Omega^2 ( \mathscr E) \otimes \Gamma (\mathscr C) \subseteq \Omega^2 (\mathscr E, T\mathscr E).
\]
In the following we will need to shift by $1$ the degree in the graded module $\overline \Omega{}^\bullet (\mathscr E, \mathscr V)$. Accordingly, we will denote by $\overline \Omega{}^\bullet (\mathscr E, \mathscr V)[1]$ the new graded module whose degree $k$ component is $\overline \Omega{}^{k+1} (\mathscr E, \mathscr V)$. We will also need the \emph{Nijenhuis-Richardson bracket} of vector valued forms \cite[Section 16]{m08} that we denote $[-,-]^{\mathrm{nr}}$. We recall that $[-,-]^{\mathrm{nr}}$ maps an $r$-form and an $s$-form into an $(r+s-1)$-form.

\begin{theorem}\label{theor:LR_infty}
A splitting (\ref{eq:split}) induces an $LR_\infty$-algebra structure on 
\[
\big(\overline \Omega{}^\bullet (\mathscr E), \overline \Omega{}^\bullet(\mathscr E, \mathscr V)[1]\big)
\]
whose structure maps are given by
\[
\begin{aligned}
\mathfrak l_1 \overline W & = \overline d\, \overline W \\
\mathfrak l_2 \big(\overline W, \overline U\big) & = - (-)^{|W|} [\overline W, \overline U]^{\mathrm{fn}} + [[R, \overline W]^{\mathrm{nr}}, \overline U]^{\mathrm{nr}} \\
\mathfrak l_3 \big(\overline W, \overline U, \overline Z\big) & = - [[[R, \overline W]^{\mathrm{nr}}, \overline U]^{\mathrm{nr}}, \overline Z]^{\mathrm{nr}}
\end{aligned}
\]
and
\[
\begin{aligned}
\mathfrak m_1 \overline \omega & = \overline d\, \overline \omega \\
\mathfrak m_2 \big(\overline W | \overline \omega\big) & = - (-)^{|W|} \mathcal L_{\overline W} \overline \omega + \iota_{[R, \overline W]^{\mathrm{nr}}} \overline \omega \\
\mathfrak m_3 \big(\overline W, \overline U| \overline \omega\big) & = - \iota_{[[R, \overline W]^{\mathrm{nr}}, \overline U]^{\mathrm{nr}}} \overline \omega
\end{aligned}
\]
for all $\overline W, \overline U, \overline Z \in \overline \Omega{}^\bullet(\mathscr E, \mathscr V)[1]$, and all $\overline \omega \in \overline \Omega{}^\bullet(\mathscr E)$, while $\mathfrak l_k = \mathfrak m_k = 0$ for $k > 3$. 
The graded Lie-Rinehart algebra structure induced in cohomology (Remark \ref{rem:LR_cohom}) agrees with that of Theorem \ref{theor:LR}.
\end{theorem}

\begin{remark}
The degree $1$ shift in the statement of Theorem \ref{theor:LR_infty} is due to our convention on $LR_\infty$-algebra structure and can be removed choosing a different convention (see \cite{vi12} for more details). 
\end{remark}

The graded $\overline \Omega{}^\bullet (\mathscr E)$-module $\Omega{}^\bullet(\mathscr E, \mathscr V)[1]$ is projective and finitely generated. Accordingly, the $LR_\infty$-algebra structure of Theorem \ref{theor:LR_infty} can be also encoded into the associated Chevalley-Eilenberg algebra. The Chevalley-Eilenberg algebra of $\big(\overline \Omega{}^\bullet (\mathscr E), \overline \Omega{}^\bullet(\mathscr E, \mathscr V)[1]\big)$ is described in the next theorem.

\begin{theorem}\label{theor:CE}
The Chevalley-Eilenberg algebra of the $LR_\infty$-algebra of Theorem \ref{theor:LR_infty} is
\[
\overline \Omega{}^\bullet (\mathscr E, \wedge^\bullet \mathscr V^\ast)
\]
with structure derivations given by
\[
\begin{aligned}
\mathfrak d_1 & = \overline d\\
\mathfrak d_2 & = d- \overline d + \iota_R \\
\mathfrak d_3 & = - \iota_R
\end{aligned}
\]
while $\mathfrak d_k = 0$ for $k > 3$. In particular the total derivation $\mathfrak D = \mathfrak d_1 + \mathfrak d_2 + \mathfrak d_3$ is just the De Rham differential $d$. The DG algebra structure induced in cohomology (Remark \ref{rem:DGA_cohom}) agrees with that of Theorem \ref{theor:DGA}.
\end{theorem}

\begin{remark}
The $LR_\infty$-algebra of Theorem \ref{theor:LR_infty} is independent of the choice of the splitting \ref{eq:split} up to \emph{$LR_\infty$-isomorphisms}. We will not provide here a notion of $LR_\infty$-morphism. We just mention that the independence of the $LR_\infty$-algebra $\big(\overline \Omega{}^\bullet (\mathscr E), \overline \Omega{}^\bullet(\mathscr E, \mathscr V)[1]\big)$ of the splitting is rigorously expressed by the fact that the total derivation $\mathfrak D = \mathfrak d_1 + \mathfrak d_2 + \cdots$ in the Chevalley-Eilenberg algebra of Theorem \ref{theor:CE} is always the same, i.e.{} the De Rham differential, regardless of the splitting we have chosen. See \cite{vi12} for more details.
\end{remark}

Theorems \ref{theor:LR_infty} and \ref{theor:CE} provide indications towards \emph{Vinogradov's conjecture} that the correct category where differential calculus on the space of solutions of a PDE $\mathscr E_0$ should be developed is the homotopy category of $DG$-modules over the horizontal De Rham algebra $(\overline \Omega{}^\bullet (\mathscr E), \overline d)$. They also suggest to add one further step to the recipe presented at the end of Section \ref{Sec:Homol}. Namely

\begin{recipe}[Addendum]\label{addendum}\quad
\begin{enumerate}
\item[(5)] Show that the cochains $\overline \Omega{}^\bullet (\mathscr E) \otimes \mathscr V \Phi$ support a homotopy algebra structure responsible for the algebraic structure on $\boldsymbol{\Phi}$.
\end{enumerate}
\end{recipe}

\subsection{the $A_\infty$-algebra of secondary differential operators}

In this final subsection we discuss \emph{secondary scalar differential operators}. In order to define them we follow the recipe presented at the end of Section \ref{Sec:Homol}. Let $(\mathscr E, \mathscr C)$ be a diffiety. The construction $\Phi$ for which we want to find a replacement on the space of solutions of $\mathscr E_0$ is ``\emph{scalar differential operators}'', i.e.{} linear differential operators from functions to functions (just DOs, for simplicity, in what follows).

For Step (1) in our recipe, we have to define \emph{transverse DOs}. Intuitively, they should be DOs taking derivatives just in the direction transverse to $\mathscr C$. Begin with the (non-commutative) $\mathbb R$-algebra $DO (\mathscr E)$ of DOs
\[
\Delta : C^\infty (\mathscr E) \to C^\infty (\mathscr E)
\]
on $\mathscr E$. As vector fields are DOs (of order $1$), sections of $\mathscr C$ span a right ideal $\Gamma (\mathscr C) \cdot DO(\mathscr E)$ in $DO (\mathscr E)$. Denote
\[
\mathscr VDO := \frac{DO(\mathscr E)}{\Gamma (\mathscr C) \cdot DO(\mathscr E)}
\] 
the quotient left $DO(\mathscr E)$-module. By definition, $\mathscr V DO$ is the space of \emph{transverse DOs}.

For Step (2) in the recipe, we have to notice that $\mathscr V DO$ is equipped with an action of the Lie algebroid $\mathscr C$. This is indeed the case: left composition with a section of $\mathscr C$ is indeed such action. Accordingly, the $\overline \Omega{}^\bullet (\mathscr E)$-module
\[
\overline \Omega{}^\bullet (\mathscr E, \mathscr V DO) := \overline \Omega{}^\bullet (\mathscr E) \otimes \mathscr V DO
\]
is a DG $(\overline \Omega{}^\bullet (\mathscr E), \overline d)$-module, whose differential we denote again by $\overline d$.

In Step (3) we define \emph{secondary differential operators} as the cohomology of $(\Omega{}^\bullet (\mathscr E, \mathscr V DO), \overline d)$:
\[
\boldsymbol{DO} := H^\bullet \big( \overline \Omega (\mathscr E, \mathscr V DO), \overline d \big).
\]

For Step (4) we have to show that $\boldsymbol{DO}$ is equipped with the appropriate algebraic structure for DOs. Actually, we can show that $\boldsymbol{DO}$ is a graded associative algebra (notice that, on the contrary, $\overline \Omega{}^\bullet (\mathscr E, \mathscr V DO)$ is not a DG algebra in general, as it does not possess a natural associative product). To do this, it is a good idea to perform simultaneously Step (5) of the recipe (see the end of the previous subsection). The relevant homotopy algebras here are \emph{$A_\infty$-algebras}.

\begin{definition}
An \emph{$A_\infty$-algebra} \cite{ls93, lm95} is a graded vector space $U$ equipped with a sequence $\mathfrak a = \{\mathfrak a_k\}_{k \in \mathbb N}$ of degree $1$ multilinear maps:
\[
\mathfrak a_k : \big(\textstyle\bigotimes^k U\big){}^\bullet \to U^{\bullet + 1}, \quad (u_1, \ldots, u_k) \mapsto \mathfrak a_k (u_1, \ldots, u_k)
\]
satisfying the following \emph{higher associativity conditions}:
\[
\sum_{r+s = k} \sum_{j=1}^{r+s}(-)^{|u_1| + \cdots + |u_j|} \mathfrak a_{s+1}\big(u_1, \ldots, u_j, \mathfrak a_{s} (u_{j+1}, \ldots, u_{j + r}), u_{j+r+1}, \ldots, u_{r+s}\big) = 0
\]
for all $k \in \mathbb N$, and all $u_1, \ldots, u_k \in U$.
\end{definition}

It follows from the above definition that $\mathfrak a_1$ is a differential, and it is a derivation with respect to $\mathfrak a_2$. Moreover $\mathfrak a_2$ is associative up to a homotopy encoded by $\mathfrak a_3$. In particular $\mathfrak a_2$ induces a honest associative product in cohomology $H^\bullet (U, \mathfrak a_1)$ (up to d\'ecalage).

\begin{theorem}\label{theor:DO_sec}
The graded vector space $\Omega{}^\bullet (\mathscr E, \mathscr V DO)[1]$ can be equipped with an $A_\infty$-algebra structure $\mathfrak a$ such that $\mathfrak a_1 = \overline d$. Moreover, the induced associative product in the cohomology $\boldsymbol{DO}$ is canonical.
\end{theorem}

The above theorem supports the interpretation of $\boldsymbol{DO}$ as \emph{secondary DOs}, i.e.{} (an appropriate replacement for) DOs on the space of solutions of $\mathscr E_0$.

We conclude the paper by sketching the proof of Theorem \ref{theor:DO_sec}. As we will see, the proof will also suggest an alternative \emph{Secondarization Recipe}. From now on we assume some familiarity with graded geometry (see, e.g., \cite{M2006}). We begin recalling a standard technique in homotopical algebra, namely the \emph{Homotopy Transfer} (see \cite{lv12, v14}). The idea is that homotopy algebras can be transferred along contraction data. A \emph{set of contraction data} is a diagram
\begin{equation}\label{eq:contr}
\xymatrix{     *{ \quad \quad \quad(A, \delta_A)\ }
\ar@(dl,ul)[]^-{h}\
\ar@<0.5ex>[r]^-{p} & *{\
(B,\delta_B)\quad\ \  \ \quad}  \ar@<0.5ex>[l]^-{j}}
\end{equation}
where 
\begin{enumerate}
\item $(A, \delta_A)$ and $(B, \delta_B)$ are cochain complexes,
\item $p, j$ are cochain maps,
\item $h : A^\bullet \to A^{\bullet -1}$ is a homotopy,
\end{enumerate}
such that
\[
pj = \operatorname{id}_B \quad \text{and} \quad jp = [\delta_A, h]
\]
and, moreover, the following \emph{side conditions} are satisfied:
\[
ph = hj = h^2 = 0.
\]
In particular $p, j$ are mutually homotopy inverse homotopy equivalences and
\[
H^\bullet(A, \delta_A) \cong H^\bullet (B, \delta_B).
\]

\begin{theorem}[Homotopy Transfer]
Let $(A, \delta_A)$ be an associative DG algebra and let (\ref{eq:contr}) be a set of contraction data over a cochain complex $(B, \delta_B)$. Then $B[1]$ can be equipped with an $A_\infty$-algebra structure $\mathfrak a$, uniquely determined by the associative product in $A$ and the contraction data, such that $\mathfrak a_1 = \delta_B$. Moreover $(A, \delta_A)$ and $(B, \mathfrak a)$ induce the same graded associative algebra structure in cohomology $H^\bullet(A, \delta_A) \cong H^\bullet (B, \delta_B)$.
\end{theorem}

There is a version of the Homotopy Transfer Theorem for $L_\infty$-algebras \cite{h10, v14} and one, slightly more involved, for $LR_\infty$-algebras (see \cite{v12} for details). Actually the $LR_\infty$-algebra of Theorem \ref{theor:LR} can be obtained via Homotopy Transfer. Namely, let $(\mathscr E, \mathscr C)$ be a diffiety. Consider the DG manifold $\mathscr C [1]$ obtained by the Lie algebroid $\mathscr C \to \mathscr E$ by shifting by one the fiber degree. Clearly $C^\infty (\mathscr C [1]) = \overline \Omega{}^\bullet (\mathscr E)$ and the cohomological vector field on $\mathscr C[1]$ is the horizontal De Rham differential $\overline d$. The pair $\big(C^\infty (\mathscr C[1]), \mathfrak X (\mathscr C[1])\big)$ is a DG Lie-Rinehart algebra. The differential in $\mathfrak X (\mathscr C[1])$ is the graded commutator of graded vector fields with the horizontal De Rham differential. Moreover, a splitting (\ref{eq:split}) uniquely determines contraction data
\begin{equation}%\label{eq:contr1}
\xymatrix{     *{ \quad \quad \quad \mathfrak X (\mathscr C[1]) \ }
\ar@(dl,ul)[]^-{h}\
\ar@<0.5ex>[r]^-{p} & *{\ \,
\overline \Omega{}^\bullet (\mathscr E, \mathscr V) \quad\ \  \ \quad}  \ar@<0.5ex>[l]^-{j}}
\end{equation}
that we can use to construct an $LR_\infty$-algebra structure on $\overline \Omega{}^\bullet (\mathscr E, \mathscr V)[1]$ from the DG Lie-Rinehart algebra structure on $(C^\infty (\mathscr C[1]), \mathfrak X (\mathscr C[1]))$ (see \cite{v12} for details). This means that the Lie-Rinehart algebra structure on secondary vector fields $(\boldsymbol{C}^\infty, \boldsymbol{\mathfrak X})$ is, equivalently, the one induced in cohomology by the DG Lie-Rinehart algebra $(C^\infty (\mathscr C[1]), \mathfrak X (\mathscr C[1]))$. A similar result holds for secondary DOs. Namely, consider the space $DO(\mathscr C[1])$ of graded DOs over the DG manifold $\mathscr C[1]$. It is an associative DG algebra whose differential is the graded commutator of graded DOs with the horizontal De Rham differential. Now, a splitting (\ref{eq:split}), together with certain additional data that we will not describe, uniquely determine contraction data
\begin{equation}%\label{eq:contr1}
\xymatrix{     *{ \quad \quad \quad \quad DO \big(\mathscr C[1]\big) \ }
\ar@(dl,ul)[]^-{h}\
\ar@<0.5ex>[r]^-{p} & *{\ \,
\overline \Omega{}^\bullet (\mathscr E, \mathscr V DO) \quad\ \  \ \quad}  \ar@<0.5ex>[l]^-{j}}
\end{equation}
that we can use to construct an $A_\infty$-algebra structure on  $\overline \Omega{}^\bullet (\mathscr E, \mathscr V DO)[1]$ from the associative DG algebra structure on $DO(\mathscr C[1])$. In particular there is a graded associative algebra structure on secondary differential operators
\[
\boldsymbol{DO} = H^\bullet \big(\overline \Omega (\mathscr E, \mathscr V DO), \overline d\big) \cong H^\bullet \big( DO (\mathscr C[1]) \big)
\]
induced by either the associative DG algebra structure on $DO (\mathscr C[1])$ or the $A_\infty$-algebra structure on $\overline \Omega{}^\bullet (\mathscr E, \mathscr V DO)$.

This discussion suggests the following alternative recipe to find the secondary analogue of a given construction $\Phi$ in differential geometry.

\begin{recipe}\quad
\begin{enumerate}
\item Apply the construction $\Phi$ to the DG manifold $\mathscr C[1]$ (when $\Phi$ is ``vector fields'', resp.{} ``DO'' we get $\mathfrak X (\mathscr C[1])$, resp.{} $DO (\mathscr C[1])$).
\item Notice that $\Phi (\mathscr C[1])$ is canonically equipped with a differential induced by the cohomological vector field $\overline d$ on $\mathscr C[1]$ (when $\Phi$ is ``vector fields'', resp.{} ``DOs'', such differential is the graded commutator of vector fields, resp.{} DOs, with $\overline d$).
\item Define $\boldsymbol{\Phi}$ as the cohomology of $\Phi (\mathscr C[1])$ (when $\Phi$ is ``vector fields'', resp. ``DOs'', then $\boldsymbol \Phi$ is $\boldsymbol{\mathfrak X}$, resp.{} $\boldsymbol{DO}$).
\item Notice that, being compatible with the differential, the algebraic structure on $\Phi (\mathscr C[1])$ induces the same algebraic structure on $\boldsymbol{\Phi}$ (when $\Phi$ is ``vector fields'', resp.{} ``DOs'', then $\boldsymbol \Phi = \boldsymbol{\mathfrak X}$, resp.{} $\boldsymbol \Phi = \boldsymbol{DO}$, is a graded Lie-Rinehart algebra, resp.{} a graded associative algebra).
\end{enumerate}
\end{recipe}

\begin{remark}
Notice that, in the alternative Secondarization Recipe just above there is no analogue of the Addendum (5) at p.~\pageref{addendum}. The reason is that the cochain complex $\Phi (\mathscr C[1])$ does possess the same (honest) algebraic structure as $\Phi$ by definition, and there is no need to work with algebraic structures up to homotopy at the level of cochains in this case. Nonetheless, by the Homotopy Transfer Theorem, cohomologies possess algebraic structures of the same type but only up to homotopy. Moreover, cohomologies with their algebraic structures up to homotopy are quasi-isomorphic to cochains with their algebraic structures. In other words, one can either choose to work with a larger space (cochains) and a simpler algebraic structure, or with a smaller space (cohomologies) and a more complicated algebraic structure (algebraic structure up to homotopy).
\end{remark}

\end{document}